\documentclass[a4paper,11pt]{article}
\usepackage{amssymb}
\usepackage{epsfig,graphics,graphics}

\newcommand{\R}{\mathbb{R}}

\newcommand{\N}{\mathbb{N}}

\newcommand{\beq}{\begin{equation} }
\newcommand{\eqq}{\end{equation} }
\newcommand{\cuad}{{\sqcap\kern-.68em\sqcup}}

\newcommand{\norm}[1]{\|#1\|}

\newtheorem{definition}{Definition}[section]
\newtheorem{teo}{Theorem}[section]

\newtheorem{proposition}{Proposition}[section]

\newtheorem{lemma}{Lemma}[section]

\newtheorem{remark}{Remark}[section]
\newcommand{\bremark}{\begin{remark} \em}
\newcommand{\eremark}{\end{remark} }

\begin{document}

\begin{center}{\bf  \Large     Semilinear fractional  elliptic equations with measures in unbounded domain }\medskip
\bigskip
{\small

{\bf Huyuan Chen\footnote{hchen@dim.uchile.cl}}
\smallskip

Departamento de Ingenier\'{\i}a  Matem\'atica and
Centro de Modelamiento Matem\'atico
 UMR2071 CNRS-UChile,
 Universidad de Chile, Chile.
 
 \medskip

{\bf Jianfu Yang\footnote{jfyang\_2000@yahoo.com} }
\smallskip

{\small Department of Mathematics, Jiangxi Normal University,\\
Nanchang, Jiangxi 330022, PR China}

}
\bigskip

\begin{abstract}
In this paper, we study the existence of  nonnegative weak solutions to (E) $ (-\Delta)^\alpha
u+h(u)=\nu $ in a general regular domain $\Omega$,
which vanish in $\R^N\setminus\Omega$, where $(-\Delta)^\alpha$
denotes the fractional Laplacian with $\alpha\in(0,1)$, $\nu$ is a nonnegative
Radon measure and $h:\R_+\to\R_+$ is a continuous nondecreasing function satisfying a subcritical integrability
condition.

Furthermore, we analyze properties of weak solution $u_k$ to $(E)$ with
 $\Omega=\R^N$, $\nu=k\delta_0$ and $h(s)=s^p$, where $k>0$, $p\in(0,\frac{N}{N-2\alpha})$ and $\delta_0$ denotes Dirac mass at the origin.
 Finally, we show
 for $p\in(0,1+\frac{2\alpha}{N}]$ that $u_k\to\infty$ in $\R^N$ as $k\to\infty$, and
 for $p\in(1+\frac{2\alpha}{N},\frac{N}{N-2\alpha})$ that $\lim_{k\to\infty}u_k(x)=c|x|^{-\frac{2\alpha}{p-1}}$ with $c>0$,
 which is a classical solution
 of  $ (-\Delta)^\alpha u+u^p=0$ in $\R^N\setminus\{0\}$.

\end{abstract}
\end{center}
  \vspace{1mm}

\vspace{2mm}
\setcounter{equation}{0}
\section{Introduction}

Let  $\Omega$ be a regular domain (not necessary bounded) of $\R^N$ ($N\geq2$),  $\alpha\in(0,1)$ and  $d\omega(x)=\frac{dx}{1+|x|^{N+2\alpha}}$.  Denote by $\mathfrak{M}^b(\Omega)$ the space of the Radon measures $\nu$ in $\Omega$
such that $\norm{\nu}_{\mathfrak{M}^b(\Omega)}:=|\nu|(\Omega)<+\infty$ and by $\mathfrak{M}^b_+(\Omega)$ the nonnegative cone.
The purpose of  this paper is to study the existence of nonnegative weak solutions to semilinear fractional elliptic equations
\begin{equation}\label{eq1.1}
 \arraycolsep=1pt
\begin{array}{lll}
 (-\Delta)^\alpha  u+h(u)=\nu\quad &{ \rm in}\quad\Omega,\\
  \phantom{   (-\Delta)^\alpha  + h(u)}
  u=0& {\rm in}\quad\R^N\setminus\Omega,
\end{array}
\end{equation}
where $h:\R_+\to\R_+$ is a continuous nondecreasing function and  $(-\Delta)^\alpha $
denotes the fractional Laplacian of exponent $\alpha$ defined by
$$(-\Delta)^\alpha  u(x)=\lim_{\epsilon\to0^+} (-\Delta)_\epsilon^\alpha u(x),$$
where for $\epsilon>0$,
$$
(-\Delta)_\epsilon^\alpha  u(x)=-\int_{\R^N}\frac{ u(z)-
u(x)}{|z-x|^{N+2\alpha}}\chi_\epsilon(|x-z|) dz
$$
and
$$\chi_\epsilon(t)=\left\{ \arraycolsep=1pt
\begin{array}{lll}
0,\quad & {\rm if}\quad t\in[0,\epsilon],\\[2mm]
1,\quad & {\rm if}\quad t>\epsilon.
\end{array}
\right.$$

In the pioneering work \cite{B12} (also see \cite{BB}),  Brezis studied the existence of weak solutions to second order elliptic problem
\begin{equation}\label{eq003}
 \arraycolsep=1pt
\begin{array}{lll}
 -\Delta  u+h(u)=\nu \quad & {\rm in}\quad\Omega,\\[2mm]
 \phantom{ -\Delta  +g(u)}
u=0  \quad & {\rm on}\quad \partial\Omega,
\end{array}
\end{equation}
where $\Omega$ is a bounded $C^2$ domain in $\R^N(N\ge3)$, $\nu$ is a bounded Radon measure in $\Omega$, and the function $h:\R\to\R$ is nondecreasing,
positive on $(0,+\infty)$ and satisfies the subcritical assumption:
$$\int_1^{+\infty}(h(s)-h(-s))s^{-2\frac{N-1}{N-2}}ds<+\infty.$$
In particular case that $0\in\Omega$, $h(s)=s^q$ and $\nu=k\delta_0$ with $k>0$,  Brezis et al showed that
(\ref{eq003}) admits  a unique weak solution $v_k$ for $1< q< N/(N-2)$, while no solution exists if $q\geq  N/(N-2)$. Later on,
V\'{e}ron in \cite{V0} proved that if $1< q< N/(N-2)$, the limit of $v_k$ is a strong singular solution of
\begin{equation}\label{local}
 \arraycolsep=1pt
\begin{array}{lll}
 -\Delta  u+u^q=0 \quad & {\rm in}\quad\Omega\setminus\{0\},\\[2mm]
 \phantom{ -\Delta  +u^q }
u=0  \quad & {\rm on}\quad \partial\Omega.
\end{array}
\end{equation}
After that,
Brezis-V\'{e}ron  in \cite{BV} found that the problem (\ref{local})
admits only the zero solution if $q\geq  N/(N-2)$. Much advances in the study
of semilinear second order elliptic equations involving measures see references \cite{BLOP,BMP1,MP,P}. 


During the last years, there has also been a renewed and increasing interest in the study of linear and nonlinear integro-differential
 operators, especially, the fractional Laplacian, motivated by various applications in physics and by important links on the theory
 of L\'{e}vy processes, refer to  \cite{CS2,CS3,CS1,CFQ,CKS,CS,CV1,CV2,FQ2,S,S1}. In a recent work, Karisen-Petitta-Ulusoy in \cite{KPU} used the duality approach to study  the fractional elliptic
equation
$$(-\Delta)^\alpha v=\mu\quad{\rm{in}}\quad \R^N,$$
where $\mu$ is a Radon measure with compact support.
More recently, Chen-V\'{e}ron in \cite{CV2}
studied the semilinear fractional elliptic problem (\ref{eq1.1}) when
$\Omega$ is an open bounded regular set in $\R^N$
 and $\nu$ is a Radon measure such that $\int_{\Omega}d^\beta d|\nu|<+\infty$ with $ \beta\in[0,\alpha]$ and $d(x)=dist(x,\partial\Omega)$.
  The existence and uniqueness of weak
solution are obtained  when the function $h$ is nondecreasing and satisfies
\begin{equation}\label{eq 2.10wx1}
\int_1^{+\infty}(h(s)-h(-s))s^{-1-k_{\alpha,\beta}}ds<+\infty,
\end{equation}
where
\begin{equation}\label{eq 2.10}
k_{\alpha,\beta}=\left\{
\arraycolsep=1pt
\begin{array}{lll}
\frac{N}{N-2\alpha},\quad &{\rm if}\quad
\beta\in[0,\frac{N-2\alpha}N\alpha],\\[2mm]
\frac{N+\alpha}{N-2\alpha+\beta},\quad &{\rm if}\quad
\beta\in(\frac{N-2\alpha}N\alpha,\alpha].
\end{array}
\right.
\end{equation}

Motivated by these results and in view of the non-local character of the fractional Laplacian we are interested in
the existence of weak solutions to problem (\ref{eq1.1}) when $\Omega$ is a general regular domain, including $\Omega=\R^N$.
Before stating our main results in this paper, we introduce the definition of weak solution to
(\ref{eq1.1}).

\begin{definition}\label{weak definition}
A function $u\in L^1(\R^N,d\omega)$ is a weak solution of (\ref{eq1.1}) if $h(u)\in L^1(\R^N,d\omega)$   and
\begin{equation}\label{weak sense}
\int_{\Omega} [u(-\Delta)^\alpha\xi+h(u)\xi]dx=\int_{\Omega} \xi d\nu,\quad \forall\xi\in \mathbb{X}_{\Omega},
\end{equation}
where  $\mathbb{X}_{\Omega}\subset C(\R^N)$ is the space of functions
$\xi$ satisfying:\smallskip

\noindent (i) the support of $\xi$ is a compact set in $\bar\Omega$;\smallskip

\noindent(ii) $(-\Delta)^\alpha\xi(x)$ exists for any $x\in\Omega$ and there exists $C>0$ such that
 $$|(-\Delta)^\alpha\xi(x)|\leq \frac{C}{1+|x|^{N+2\alpha}},\quad \forall x\in\Omega;$$

\noindent(iii) there exist $\varphi\in L^1(\Omega,\rho^\alpha dx)$
and $\epsilon_0>0$ such that $|(-\Delta)_\epsilon^\alpha\xi|\le
\varphi$ a.e. in $\Omega$, for all
$\epsilon\in(0,\epsilon_0]$, here $\rho(x)=\min\{1,dist(x,\partial\Omega)\}$ if $\Omega\not=\R^N$ or $\rho\equiv1$ if $\Omega=\R^N$.\smallskip
\end{definition}

We notice that $\mathbb{X}_{\Omega}$ coincides with the test function space  $\mathbb{X}_{\alpha}$ if $\Omega$ is bounded, see \cite[Definition 1.1]{CV2}.
Moreover, the test function space $\mathbb{X}_{\Omega}$ is used as $C^{1,1}_0(\Omega)$ if $\Omega$ is bounded and $\alpha=1$, see \cite{V}.
We denote by  $G_\Omega$ the Green kernel of $(-\Delta)^\alpha$ in
$\Omega\times\Omega $ and by $\mathbb{G}_\Omega[\cdot]$ the Green operator
defined as
$$
\mathbb{G}_\Omega[\nu](x)=\int_{\Omega}G_\Omega(x,y) d\nu(y),\quad\forall \nu\in
\mathfrak{M}^b(\Omega).
$$

Now we are ready to state our first theorem on the existence of weak solutions for problem (\ref{eq1.1}).
\begin{teo}\label{teo 1}
Assume  that  $\alpha\in(0,1)$, $\Omega$ is a regular domain in $\R^N (N\ge2)$ and  $h:\R_+\to\R_+$
is a continuous  nondecreasing function  satisfying
\begin{equation}\label{1.4}
 \int_1^{+\infty} h(s) s^{-1-\frac{N}{N-2\alpha}}ds <+\infty.
\end{equation}
Then for any $\nu\in\mathfrak{M}^b_+(\Omega)$, problem
(\ref{eq1.1}) admits a weak solution $u_\nu$ such that
\begin{equation}\label{1.5}
0\le u_\nu\le \mathbb{G}_\Omega[\nu]\quad \rm{a.e.\ in}\
\Omega.
\end{equation}
\end{teo}

If $\nu$ is a nonnegative bounded Radon measure, (\ref{eq 2.10wx1})  with $\beta=0$ and
(\ref{1.4}) have the same critical value $\frac{N}{N-2\alpha}$.

In the case that $\Omega$ is bounded, the authors of \cite{CV2} took a sequence of $C^1$ functions $\{\nu_n\}$ converging to $\nu$ in the weak star sense, then they considered the solutions $u_n$ of (\ref{eq1.1}) replacing $\nu$ by $\nu_n$. By the compact imbedding theorem, they showed that the limit of $\{u_n\}$ exists, up to subsequence. While for the case that $\Omega$ is unbounded, the difficulty is that Sobolev imbedding may not be compact. To overcome the difficulty, we truncate the measure $\nu$ by $\nu\chi_{B_R(0)}$ and use the increasing monotonicity of corresponding solutions
$\{u_R\}$ of solutions to (\ref{eq1.1}) in related bounded domains. Taking the limit as $R\to\infty$, we achieve the desired weak solution.

The second purpose in this paper is to study properties of weak solution to problem (\ref{eq1.1}) when $\Omega=\R^N$, $h(u)=u^p$ and $\nu=k\delta_0$, that is,
\begin{equation}\label{eq1.3}
 \arraycolsep=1pt
\begin{array}{lll}
 (-\Delta)^\alpha  u+u^p=k\delta_0\quad   {\rm in}\quad \R^N,\\[2mm]
  \phantom{    }
  \lim_{|x|\to+\infty}u(x)=0,
\end{array}
\end{equation}
where $p\in(0,\frac{N}{N-2\alpha})$, $k>0$ and $\delta_0$ denotes the Dirac mass at the origin.

\begin{teo}\label{cor 1}
Assume that $\alpha\in(0,1)$ and  $p\in(0,\frac{N}{N-2\alpha})$.
Then for any $k>0$, problem
(\ref{eq1.3}) admits a unique weak solution $u_k$ such that
\begin{equation}\label{1.5}
\lim_{x\to0}u_k(x)|x|^{N-2\alpha}=c_1k,
\end{equation}
where $c_1>0$.
 Moreover,
\begin{itemize}
\item[]
\begin{enumerate}\item[$(i)$]
$\{u_k\}_{k\in(0,\infty)}$  are classical solutions of
\begin{equation}\label{eq1.01}
 (-\Delta)^\alpha  u+u^p=0\quad { \rm in}\quad\R^N\setminus\{0\};
\end{equation}
\end{enumerate}
 \begin{enumerate}\item[$(ii)$]
the mapping: $k\mapsto u_k$ is increasing.
\end{enumerate}
\end{itemize}
\end{teo}

We consider the asymptotic behavior of $u_1$ at $\infty$  when $p\in(1,\frac{N}{N-2\alpha})$.

\begin{teo}\label{th 1}
Assume that $\alpha\in(0,1)$ and $u_1$ is the solution of (\ref{eq1.3}) with $k=1$.
Then there exist $c_2>1$ and $R>2$ such that for $|x|\ge R$,
\begin{itemize}
\item[]
\begin{enumerate}\item[$(i)$]
if $p\in(1,1+\frac{2\alpha}{N})$,
\begin{equation}\label{6.1}
\frac{1}{c_2}|x|^{-\frac{N+2\alpha}{p}}\le u_1(x)\le c_2|x|^{-\frac{N+2\alpha}{p}};
\end{equation}
\end{enumerate}
 \begin{enumerate}\item[$(ii)$]
if $p=1+\frac{2\alpha}{N}$,
\begin{equation}\label{6.3}
  \frac{1}{c_2}|x|^{-N}\log^{\frac{N}{2\alpha}}(|x|)\le u_1(x)\le c_2|x|^{-N}\log^{\frac{N}{2\alpha}}(|x|);
\end{equation}
\end{enumerate}
 \begin{enumerate}\item[$(iii)$]
if $p\in(1+\frac{2\alpha}{N}, \frac{N}{N-2\alpha})$,
\begin{equation}\label{6.2}
\frac{1}{c_2}|x|^{-\frac{2\alpha}{p-1}}\le u_1(x)\le c_2|x|^{-\frac{2\alpha}{p-1}}.
\end{equation}
\end{enumerate}
\end{itemize}

\end{teo}

According to Theorem \ref{th 1}, we know that the decaying power of $u_1$ shifts at the point $p=1+\frac{2\alpha}{N}$;
while for $\alpha=1$ and $p\in(1,\frac{N}{N-2})$, the weak solution of (\ref{eq1.3}) decays   as $|x|^{-\frac{2}{p-1}}$.

From now on, we denote that $u_k$ is the weak solution of (\ref{eq1.3}).
Since the mapping: $k\mapsto u_k$ is increasing by Theorem \ref{cor 1}
and then  we can denote that
\begin{equation}\label{definition infty}
  u_{\infty}(x)=\lim_{k\to\infty}u_k(x), \quad  x\in\R^N.
\end{equation}
Here we note that $u_{\infty}(x)\in \R_+\cup\{+\infty\}$ for any $x\in\R^N$.
Now we state properties of $u_\infty$.

\begin{teo}\label{teo 2}
Assume that  $\alpha\in(0,1)$,   $p\in(0,\frac{N}{N-2\alpha})$  and $u_\infty$ is given by (\ref{definition infty}).
Then
\begin{itemize}
\item[]
\begin{enumerate}\item[$(i)$]
if $p\in(0,1+\frac{2\alpha}N]$, then $ u_\infty(x)=\infty, \  \forall x\in \R^N;$
\end{enumerate}
 \begin{enumerate}\item[$(ii)$]
if $p\in(1+\frac{2\alpha}N, \frac{N}{N-2\alpha})$, then $u_\infty$ is a classical solution of
(\ref{eq1.01}) and there exists $c_3>0$ such that
$$u_\infty(x)=c_3|x|^{-\frac{2\alpha}{p-1}},\quad \forall x\in\R^N\setminus\{0\}.$$

\end{enumerate}
\end{itemize}
\end{teo}

In the proof of Theorem \ref{teo 2}, we make use of the self-similar property of $u_\infty$.

Analogue results of Theorem \ref{teo 2} in bounded domain $\Omega$ were obtained \cite{CV2,CV3}. Precisely,
they showed that there exists a unique weak solution $u_{k,\Omega}$ to
  semilinear fractional elliptic problem
\begin{equation}\label{eq1.001}
 \arraycolsep=1pt
\begin{array}{lll}
 (-\Delta)^\alpha  u+u^p=k\delta_0\quad & {\rm in}\quad\Omega,\\[2mm]
 \phantom{   (-\Delta)^\alpha  + u^p}
u=0\quad & {\rm in}\quad \Omega^c,
\end{array}
\end{equation}
where $k>0$, $0\in \Omega$ and $p\in(0,\frac{N}{N-2\alpha})$.
Moreover,\\
   $(i)$ the mapping $k\mapsto u_{k,\Omega}$ is increasing;\\
 $(ii)$ for $p\in(0,\min\{1+\frac{2\alpha}N, \frac{N}{2\alpha}\})$,
 $u_{\infty,\Omega}=\infty$ in $\Omega,$
 where
 $$u_{\infty,\Omega}=\lim_{k\to\infty}u_{k,\Omega}\quad{\rm in}\ \ \R^N;$$
 $(iii)$ for $p\in(1+\frac{2\alpha}N, \frac{N}{N-2\alpha})$, $u_{\infty,\Omega}$ is a classical solution of
\begin{equation}\label{eq1.02}
\arraycolsep=1pt
\begin{array}{lll}
 (-\Delta)^\alpha  u+u^p=0\quad  &{\rm in}\quad\Omega\setminus\{0\},\\[2mm]
 \phantom{   (-\Delta)^\alpha  + u^p}
 u=0\quad & {\rm in}\quad \Omega^c.
 \end{array}
\end{equation}

Finally, we discuss properties of weak solution $u_{k,\Omega}$ of (\ref{eq1.001})
when $\Omega$ is an unbounded regular domain including the origin. The result is stated as follows.

\begin{teo}\label{teo 3}
Assume that $\alpha\in(0,1)$, $\Omega$ is an unbounded regular domain of $\R^N(N\ge2)$ including the origin,
 $p\in(0,\frac{N}{N-2\alpha})$ and
$u_k$ is given by Theorem \ref{cor 1}.
Then

 $(i)$ for any $k>0$, (\ref{eq1.001}) admits a unique weak solution $u_{k,\Omega}$
such that
$$u_k-m_{k,\Omega}\le  u_{k,\Omega}\le u_k\quad{\rm in}\quad \Omega$$
and the mapping $k\mapsto u_{k,\Omega}$ is increasing,
where $m_{k,\Omega}=\sup_{x\in\Omega^c} u_k(x)$;

 $(ii)$ for $p\in(0,\min\{1+\frac{2\alpha}N, \frac{N}{2\alpha}\})$, $u_{\infty,\Omega}=\infty$ in $\Omega,$
 where $$u_{\infty,\Omega}=\lim_{k\to\infty}u_{k,\Omega}\quad{\rm in}\ \ \R^N;$$

$(iii)$  for $p\in(1+\frac{2\alpha}N, \frac{N}{N-2\alpha})$, $u_{\infty,\Omega}$ is a classical solution of
(\ref{eq1.02}) such that
$$u_\infty-m_{\infty,\Omega}\le  u_{\infty,\Omega}\le u_\infty\quad{\rm in}\quad \Omega,$$
where $u_\infty$ is defined by (\ref{definition infty}) and $m_{\infty,\Omega}=\sup_{x\in\Omega^c} u_\infty(x).$
\end{teo}

The paper is organized as follows. In Section 2 we list some
properties of Marcinkiewicz spaces and establish the inequality
\begin{equation}\label{1.3'}
\|\mathbb{G}_{\R^N}[\nu]\|_{M^{\frac{N}{N-2\alpha}}(\R^N, d\omega)}\le
c_5\|\nu\|_{\mathfrak{M}^b(\R^N)},
\end{equation}
which is used to obtain that $h(\mathbb{G}_{\R^N}[\nu])\in L^1(\R^N, d\omega)$.
 In Section  3, we prove Theorem \ref{teo 1}.
The proofs of Theorem \ref{cor 1} and Theorem \ref{th 1} are addressed in  Section 4. 
Finally, we give the proofs of  Theorem \ref{teo 2} and Theorem \ref{teo 3} in Section 5.

\setcounter{equation}{0}
\section{Preliminary}

The purpose of this section is to introduce some preliminaries and prove
Marcinkiewicz type estimate.

\subsection{Marcinkiewicz type estimate}

In this subsection,
we recall the definition of Marcinkiewicz space and  prove
Marcinkiewicz type estimate.

\begin{definition}
Let $\Theta\subset \R^N$ be a domain and $\mu$ be a positive
Borel measure in $\Theta$. For $\kappa>1$,
$\kappa'=\kappa/(\kappa-1)$ and $u\in L^1_{loc}(\Theta,d\mu)$, we
set
\begin{equation}\label{mod M}
\|u\|_{M^\kappa(\Theta,d\mu)}=\inf\left\{c\in[0,\infty]:\int_E|u|d\mu\le
c\left(\int_Ed\mu\right)^{\frac1{\kappa'}},\ \forall E\subset \Theta,\,E\
\rm{Borel}\right\}
\end{equation}
and
\begin{equation}\label{spa M}
M^\kappa(\Theta,d\mu)=\{u\in
L_{loc}^1(\Theta,d\mu):\|u\|_{M^\kappa(\Theta,d\mu)}<\infty\}.
\end{equation}
\end{definition}

$M^\kappa(\Theta,d\mu)$ is called the Marcinkiewicz space of
exponent $\kappa$, or weak $L^\kappa$-space and
$\|.\|_{M^\kappa(\Theta,d\mu)}$ is a quasi-norm.

\begin{proposition}\label{pr 1} \cite{BBC}
Assume that $1\le q< \kappa<\infty$ and $u\in L^1_{loc}(\Theta,d\mu)$.
Then there exists  $c_4>0$ dependent of $q,\kappa$ such that
$$\int_E |u|^q d\mu\le c_4\|u\|_{M^\kappa(\Theta,d\mu)}\left(\int_E d\mu\right)^{1-q/\kappa},$$
for any Borel set $E$ of $\Theta$.
\end{proposition}

Now we are ready to state Marcinkiewicz type estimate as follows.

\begin{proposition}\label{general}
Let  $\nu\in\mathfrak{M}^b(\R^N)$, then there exists $c_5>0$ such that
\begin{equation}\label{annex 0}
\|\mathbb{G}_{\R^N}[|\nu|]\|_{M^{p^*_\alpha}(\R^N,d\omega)}\le c_5\|\nu\|_{\mathfrak{M}^b(\R^N)},
\end{equation}
where $d\omega(x)=\frac{dx}{1+|x|^{N+2\alpha}}$ and $p^*_\alpha=\frac{N}{N-2\alpha}$.

\end{proposition}
{\it Proof}.  For $\lambda>0$ and $y\in \R^N$, we set
$$A_\lambda(y)=\left\{x\in \R^N\setminus\{y\}: G_{\R^N}(x,y)>\lambda \right\},\quad
 m_\lambda(y)=\int_{A_\lambda(y)}d\omega.$$
We observe that there exists a positive constant $c_{N,\alpha}$ such that
$$
G_{\R^N}(x,y)= \frac{c_{N,\alpha}}{|x-y|^{N-2\alpha}},\quad (x,y)\in
\R^N\times\R^N,\ x\neq y,
$$
which implies that
\begin{equation}\label{gra 0003}
A_\lambda(y)\subset\left\{x\in \R^N: |x-y|\le c_{N,\alpha}\lambda^{-\frac1{N-2\alpha}}\right\}.
\end{equation}
As a consequence,
\begin{equation}\label{gra 0003wx2}
 m_\lambda(y)=\int_{A_\lambda(y)}\frac{dx}{1+|x|^{N+2\alpha}}\le|A_\lambda(y)| \le c_6\lambda^{-p_\alpha^*},
 \end{equation}
where $c_6>0$ independent of $y$ and $\lambda$ and $p_\alpha^*=\frac{N}{N-2\alpha}$.

Let $E\subset \R^N$  be a Borel set and $\lambda>0$, then
\begin{eqnarray*}
\int_E G_{\R^N}(x,y) d\omega(x)\le
\int_{A_\lambda(y)} G_{\R^N}(x,y) d\omega(x)+\lambda\int_Ed\omega.
\end{eqnarray*}
By (\ref{gra 0003wx2}), we have that
\begin{eqnarray*}
\int_{A_\lambda(y)}  G_{\R^N}(x,y)d\omega(x)=\lambda m_\lambda(y)+ \int_{\lambda}^\infty m_s(y)ds\le c_7\lambda^{1-p_\alpha^*},
\end{eqnarray*}
for some $c_7>0$, then it results that
\begin{eqnarray*}
\int_E  G_{\R^N}(x,y)d\omega(x)\le c_7 \lambda^{1-p_\alpha^*}+\lambda \int_Ed\omega.
\end{eqnarray*}
Choosing $\lambda=(\int_Ed\omega)^{-\frac1{p_\alpha^*}}$, we obtain
\begin{eqnarray*}
\int_E  G_{\R^N}(x,y)d\omega(x)\le (c_7+1)\left(\int_E d\omega\right)^{\frac{p_\alpha^*-1}{p_\alpha^*}},\quad \forall y\in \R^N.
\end{eqnarray*}
Therefore,
\begin{equation}\label{none}\begin{array}{lll}\displaystyle
\displaystyle \int_E \mathbb{G}_{\R^N}[|\nu|](x)d\omega(x)=\int_{\R^N}\int_E
 G_{\R^N}(x,y)  d\omega(x) d|\nu(y)|
\\[3mm]\phantom{---------}
\displaystyle\le (c_7+1)\int_{\R^N} d|\nu(y)|\left(\int_E d\omega\right)^{\frac{p^*_\alpha-1}{p^*_\alpha}}
\\[3mm]\phantom{---------}
\displaystyle\le (c_7+1)\|\nu\|_{\mathfrak{M}^b(\R^N)} \left(\int_E
d\omega\right)^{\frac{p^*_\alpha-1}{p^*_\alpha}}.
\end{array}
\end{equation}
As a consequence,
\begin{eqnarray*}
\| \mathbb{G}_{\R^N}[|\nu|]\|_{M^{p^*_\alpha}(\R^N,d\omega)}\le
c_5\|\nu\|_{\mathfrak{M}^b(\R^N)},
\end{eqnarray*}
which ends the proof.\qquad$\Box$

Now we use Marcinkiewicz type estimate  to prove the following lemma, which is the key-stone in the proof of Theorem \ref{teo 1}.

\begin{lemma}\label{lm 2.2}
Assume that  $\nu\in \mathfrak{M}_+^b(\R^N)$ and $h:\R_+\to\R_+$
is a continuous nondecreasing function satisfying
(\ref{1.4}).
Then
$$\mathbb{G}_{\R^N}[\nu],\ h(\mathbb{G}_{\R^N}[\nu]) \in L^1(\R^N, d\omega).$$

\end{lemma}
{\bf Proof.} On the one hand,  using Fubini's lemma, we have that
$$
 \arraycolsep=1pt
\begin{array}{lll}
  \norm{\mathbb{G}_{\R^N}[\nu]}_{L^1(\R^N,d\omega)}=c_{N,\alpha}\int_{\R^N}\int_{\R^N}\frac1{1+|x|^{N+2\alpha}}\frac{1}{|x-y|^{N-2\alpha}}d\nu(y)dx \\[3mm]\phantom{--------\ }
= c_{N,\alpha}\int_{\R^N}\int_{\R^N}\frac1{1+|x|^{N+2\alpha}}\frac{1}{|x-y|^{N-2\alpha}}dx d\nu(y)
 \\[3mm]\phantom{--------\ }
\le c_{N,\alpha}\int_{\R^N}[\int_{B_1(y)}\frac{1}{|x-y|^{N-2\alpha}}dx+\int_{B_1^c(y)}\frac1{1+|x|^{N+2\alpha}}dx ]d\nu(y)
 \\[3mm]\phantom{--------\ }
 <+\infty,
\end{array}
$$
that is, $\mathbb{G}_{\R^N}[\nu]\in L^1(\R^N, d\omega).$

On the other hand,   let $S_\lambda=\{x\in B_R(0):\mathbb{G}_{\R^N}[\nu](x)>\lambda\}$  and
$g(\lambda)=\int_{S_\lambda} d\omega$, where $\lambda\ge1$. We observe that
\begin{equation}\label{2.4}
\displaystyle\begin{array}{lll}
\displaystyle
\int_{\R^N}h(\mathbb{G}_{\R^N}[\nu])d\omega=\int_{S^c_{\lambda}}h(\mathbb{G}_{\R^N}[\nu])d\omega+\int_{S_{\lambda}}
h(\mathbb{G}_{\R^N}[\nu])d\omega
\\[4mm]\phantom{-\int_{E}|g(u_{n_k})|\rho^{\beta}dx}
\leq h(\lambda)\int_{\R^N}d\omega+\int_{S_{\lambda}}h(\mathbb{G}_{\R^N}[\nu])d\omega
\\[4mm]\phantom{-\int_{E}|g(u_{n_k})|\rho^{\beta}dx}
= h(\lambda)\int_{\R^N} d\omega+h(\lambda)g(\lambda)+\int_\lambda^{\infty}g(s)dh(s).
\end{array}
\end{equation}
Since
$$\int_{\lambda}^\infty g(s) dh(s)=\lim_{T\to\infty}\int_{\lambda}^T g(s)dh(s)$$
and $\mathbb{G}_{\R^N}[\nu]\in M^{p^*_{\alpha}}(\R^N,d\omega)$,  it derives from Proposition \ref{general} and Proposition \ref{pr 1} with $q=1$, $\kappa=p^*_\alpha$, $E=S_\lambda$
and $d\mu=d\omega$ that
$g(s)\leq c_{8}s^{-p^*_{\alpha}}$ and for $T>\lambda$,
$$\displaystyle\begin{array}{lll}
h(\lambda)g(\lambda)+\int_\lambda^Tg(s)dh(s)\le c_{8}\lambda^{-p^*_{\alpha}} h(\lambda)+ c_{8}\int_{\lambda}^T s^{-p^*_{\alpha}}dh(s)
\\[4mm]\phantom{-----------\ \ }\displaystyle
\leq c_8T^{-p^*_\alpha}h(T)+c_{8}p^*_{\alpha}\int_{\lambda}^T
s^{-1-p^*_{\alpha}}h(s)ds,
\end{array}$$
where $c_8>0$. By the nondecreasing monotonicity of $h$, we have that
\begin{eqnarray*}
T^{-p^*_{\alpha}}h(T) &=& 2^{1+p_\alpha^*}h(T)(2T)^{-1-p^{*}_\alpha}\int_T^{2T}dt \\
   &\le& 2^{1+p_\alpha^*}\int_T^{2T}h(t)t^{-1-p^*_\alpha}dt,
\end{eqnarray*}
then it infers by  (\ref{1.4}) that
that $$\lim_{T\to\infty}T^{-p^*_{\alpha}}h(T)= 0.$$
Therefore, by  (\ref{1.4}) and we take $\lambda=1$,
$$\int_{\R^N}h(\mathbb{G}_{\R^N}[\nu])d\omega\leq h(1)\int_{\R^N} d\omega+c_{8}p^*_{\alpha}\int_{1}^\infty s^{-1-p^*_{\alpha}}h(s)ds<+\infty,
$$
i.e. $h(\mathbb{G}_{\R^N}[\nu]) \in L^1(\R^N, d\omega).$
We complete the proof.\qquad$\Box$

\subsection{Basic results}

This subsection is devoted to present some basic results and Comparison Principle, which are key tools in the analysis.
We start it by recalling the existence of weak solution of (\ref{eq1.1}) when $\Omega$ is a bounded $C^2$ domain.

\begin{proposition}\cite[Theorem 1.1]{CV2}\label{pr 2.1}
Assume that $\mathcal{O}$ is a bounded $C^2$ domain in $\R^N$,  $\mu\in\mathfrak{M}^b_+(\mathcal{O})$ and  $h:\R_+\to\R_+$
is a continuous nondecreasing function satisfying (\ref{1.4}). Then problem
\begin{equation}\label{2.1}
 \arraycolsep=1pt
\begin{array}{lll}
 (-\Delta)^\alpha  u+h(u)=\mu\quad & {\rm in}\quad \mathcal{O},\\
  \phantom{   (-\Delta)^\alpha  + h(u)}
  u=0& {\rm in}\quad\R^N\setminus\mathcal{O}
\end{array}
\end{equation}
 admits a unique weak solution $v_\mu$ such that
\begin{equation}\label{2.2}
0\le v_\mu\le \mathbb{G}_{\mathcal{O}}[\mu]\quad{\rm a.e.\ in}\ \mathcal{O}.
\end{equation}
Moreover, the mapping $\mu\mapsto v_\mu$ is increasing.

\end{proposition}

Next we recall the Comparison Principle from \cite{CFQ}.

\begin{lemma}\label{teo CP}\cite[Theorem 2.3]{CFQ}
Suppose that $O$ is a bounded domain of $\R^N$, $p>0$, the functions
 $u_1$, $u_2$ are continuous in $\bar O$ and    satisfy
$$(-\Delta)^\alpha  u_1+|u_1|^{p-1}u_1\ge0 \ {\rm in}\  O\quad {\rm and} \quad (-\Delta)^\alpha  u_2+|u_2|^{p-1}u_2\le0 \ {\rm in}\  O.$$
If $u_1\ge u_2$ a.e. in $O^c$, then   $u_1\ge u_2$ in $O$.
\end{lemma}

By the Comparison Principle, we have the following result:

\begin{lemma}\label{lm 2.3}
 Assume that $f\in C^1(\R^N)$ is a nonnegative function, $h$ is a continuous and nondecreasing function  and $\mathcal{O}_1,\mathcal{O}_2$ are bounded $C^2$ domain such that $\mathcal{O}_1\subset \mathcal{O}_2$.  Let $w_1$ and $w_2$ be  the solutions of
(\ref{2.1}) with $\mu=f$ in  $\mathcal{O}=\mathcal{O}_1$ and   $\mu=f$ in $\mathcal{O}=\mathcal{O}_2$, respectively.
Then $$w_1\le w_2\quad {\rm in}\quad \R^N.$$

\end{lemma}
{\bf Proof.} Since $\mathcal{O}_1\subset \mathcal{O}_2$ and $f\ge0$, it follows by Lemma \ref{teo CP} that
 $w_2\ge 0$ in $\mathcal{O}_2$. Suppose on the contrary that
$$\min_{x\in \R^N}(w_2-w_1)(x)<0,$$
there would exist $x_0\in \mathcal{O}_1$ such that
$$(w_2-w_1)(x_0)=\min_{x\in \R^N}(w_2-w_1)(x).$$
 Then $$
 \arraycolsep=1pt
\begin{array}{lll}
(-\Delta)^\alpha(w_2-w_1)(x_0)=-\lim_{\epsilon\to0^+}\int_{\R^N\setminus B_\epsilon(x_0)}\frac{(w_2-w_1)(z)-(w_2-w_1)(x_0)}{|z-x_0|^{N+2\alpha}}dz<0
\end{array}
$$
and $h(w_2(x_0))\leq h(w_1(x_0))$,
 which  implies a contradiction since $w_1$ and $w_2$ satisfy
 $(-\Delta)^\alpha u(x_0)+h(u(x_0))=f(x_0).$
 The proof is completed.\qquad$\Box$

\begin{lemma}\label{lm 2---1}
 Suppose that $p>0$, ${\mathcal{O}}$ is a bounded $C^2$ domain in $\R^N$,  $g\in L^1({\mathcal{O}}^c,d\omega)$  is  $C^2$ in $\{z\in{\mathcal{O}}^c, {\rm dist}(z,\partial{\mathcal{O}})\leq \delta \}$ with $\delta>0$.
Then there exists a unique classical solution $u$ of
\begin{equation}\label{3.2.4}
\left\{ \arraycolsep=1pt
\begin{array}{lll}
 (-\Delta)^{\alpha} u(x)+|u|^{p-1}u(x)=0,\ \ \ \ &
x\in{\mathcal{O}},\\[2mm]
u(x)=g(x),\ &x\in{\mathcal{O}}^c.
\end{array}
\right.
\end{equation}
\end{lemma}
{\bf Proof.} For the existence of classical solutions, we refer to Theorem 2.5 in \cite{CFQ}.
The uniqueness follows by Lemma \ref{teo CP}.\qquad$\Box$

\setcounter{equation}{0}
\section{Existence of weak solutions}

In this section, we show the existence of solutions of problem (\ref{eq1.1}), that is, we will prove Theorem \ref{teo 1}. We first give an auxiliary lemma as follows.

\begin{lemma}\label{lm 2.1}
Assume that $\mathcal{O} $ is a  bounded $C^2$ domain in $\R^N$ and  $\eta\in C(\R^N)$ with support in $\bar \mathcal{O}$.
Then there exists $c_9>0$ such that
\begin{equation}\label{2.3}
 |(-\Delta)^\alpha\eta(x)|\le  \frac {c_9\|\eta\|_{L^\infty(\mathcal{O})}}{1+|x|^{N+2\alpha}},\quad x\in \R^N\setminus\mathcal{O}_d,
\end{equation}
where $\mathcal{O}_d=\{x\in \R^N: dist(\mathcal{O},x)\le d\}$.
\end{lemma}
{\bf Proof.} For $x\in \mathcal{O}_d^c$ and $y\in \mathcal{O}$, there exists $c_{10}>1$ such that
$$c_{10}^{-1}(1+|x|^{N+2\alpha})\le |x-y|^{N+2\alpha}\le c_{10}( 1+|x|^{N+2\alpha}).$$
By the fact that
$$
 (-\Delta)^\alpha\eta(x) = -\int_{\mathcal{O}}\frac{\eta(y)}{|x-y|^{N+2\alpha}}dy,\quad x\in \mathcal{O}_d^c,
$$
we assert (\ref{2.3}) holds, which ends the proof. \qquad$\Box$\medskip

Now we are in the position to prove Theorem \ref{teo 1}.
\smallskip

\noindent{\bf Proof of Theorem \ref{teo 1}.} Let $\{\mathcal{O}_n\}_{n\in\N}$ be a sequence of $C^2$ domains in $\R^N$ such that
$$\Omega\cap B_{n}(0)\subset\{\mathcal{O}_n\}\subset \Omega\cap B_{n+1}(0).$$
For $\nu\in\mathfrak{M}^b_+(\Omega)$ and $n\in\N$, we denote $\nu_n=\nu\chi_n$,
where $\chi_n$ is the characteristic function in $\mathcal{O}_n$.
By Proposition \ref{pr 2.1},  problem (\ref{2.1}) with $\mathcal{O}=\mathcal{O}_n$ and $\mu=\nu_{n}$
admits a unique weak solution, denoting it by $v_n$. We divide the proof in following steps.

{\it Step 1. We claim that $v_{n}\le v_{n+1}$ a.e. in $\R^N$.}
In fact, let $\tilde v_{n+1}$ be the solution of (\ref{2.1}) with $\mathcal{O}=\mathcal{O}_{n+1}$ and $\mu=\nu_{n}$.
By Proposition \ref{pr 2.1},
\begin{equation}\label{3.1}
\tilde v_{n+1}\le v_{n+1}\quad{\rm a.e.\ in}\ \R^N.
\end{equation}

Choosing a sequence nonnegative functions $\{f_{n,m}\}_{m\in\N}\subset C^1_0(\mathcal{O}_{n})$
such that $f_{n,m}\rightharpoonup \nu_{n}$ as $m\to\infty$ in the distribution sense, we make zero extension of $f_{n,m}$ into $C_0^1(\mathcal{O}_{n+1})$
and denote the extension by $\tilde f_{n,m}$. Let $v_{n,m}$ and $\tilde v_{n+1,m}$ be solutions
of (\ref{2.1}) with $\mu=f_{n,m}$, $\mathcal{O}=\mathcal{O}_n$ and  $\mu=\tilde f_{n,m}$, $\mathcal{O}=\mathcal{O}_{n+1}$, respectively. Lemma \ref{lm 2.3} implies that
$$
v_{n,m}\le \tilde v_{n+1,m}\quad{\rm in}\ \R^N.
$$
Together with the facts that
 $v_{n,m}\to v_n$ a.e. in $\R^N$ and $\tilde v_{n+1,m}\to \tilde v_{n+1}$ a.e. in $\R^N$ as $m\to\infty$, we obtain that
\begin{equation}\label{3.2}
v_{n}\le \tilde v_{n+1}\quad{\rm a.e.\ in}\ \ \R^N.
\end{equation}
It follows by (\ref{3.1}) and (\ref{3.2}) that for any $n\in\N$,
\begin{equation}\label{3.3}
 v_{n}\le v_{n+1}\quad {\rm a.e.\ in}\ \ \R^N.
\end{equation}

{\it Step 2. Uniform bounds of $\{v_{n}\}$.} We deduce by (\ref{2.2}) that
\begin{equation}\label{3.4}
0\le v_{n}\le \mathbb{G}_{\mathcal{O}_{n}}[\nu_{n-1}]\quad {\rm a.e.\ in}\ \ \R^N.
\end{equation}
Observing that for any $n\in\N$,
$$G_{\mathcal{O}_{n}}(x,y)\le G_{\Omega}(x,y)\le  G_{\R^N}(x,y),\quad {\rm for\ any}\ x,y\in\R^N,x\not=y$$
and $\nu_{n-1}\le \nu$, we have that
\begin{eqnarray*}
 \mathbb{G}_{\mathcal{O}_{n}}[\nu_{n-1}](x)
   \le \int_{\Omega}G_{\Omega}(x,y)d\nu(y) =\mathbb{G}_{\Omega}[\nu](x)\le \mathbb{G}_{\R^N}[\nu](x),\quad x\in\R^N,
\end{eqnarray*}
where we make zero extension of $\nu$ such that $\nu\in\mathfrak{M}^b(\R^N)$.
Therefore, by (\ref{3.4}) we obtain that  for any $n\in\N$,
\begin{equation}\label{3.5}
v_n\le \mathbb{G}_{\Omega}[\nu]\le \mathbb{G}_{\R^N}[\nu]\quad {\rm a.e.\ in}\ \ \R^N.
\end{equation}

{\it Step 3. Existence of weak solution.}
By (\ref{3.3}) and (\ref{3.5}), we see that the limit of  $\{v_n\}$ exists, denoted it by $u_\nu$. Hence,
\begin{equation}\label{2.1.1}
0\le u_\nu\le \mathbb{G}_{\Omega}[\nu]\le \mathbb{G}_{\R^N}[\nu]\quad {\rm a.e.\ in}\ \ \R^N.
\end{equation}
It follows by Lemma \ref{lm 2.2} that $\mathbb{G}_{\R^N}[\nu]\in L^1(\R^N, d\omega)$ and then $u_\nu\in L^1(\R^N, d\omega)$.
Thus, $v_n\to u_\nu$  in $L^1(\R^N, d\omega)$ as $n\to\infty$.
Moreover,
\begin{itemize}
\item[]
\begin{enumerate}\item[$(i)$]
$\{h(v_n)\}_{n\in\N}$ is an increasing sequence of functions and $h(v_n)\to h(u_\nu)$ a.e. in $\R^N$;
\end{enumerate}
 \begin{enumerate}\item[$(ii)$]
 it implies by $\Omega^c\subset\mathcal{O}_n^c$ and $v_n=0$ in $\mathcal{O}_n^c$ that $u_\nu=0$ in $\Omega^c$.
\end{enumerate}
\end{itemize}

For $\xi\in \mathbb{X}_{\Omega}$, there exists $N_0>0$ such that for any $n\ge N_0$,
 $${\rm supp}(\xi)\subset \bar\mathcal{O}_n,$$
which implies that
$\xi\in \mathbb{X}_{\mathcal{O}_n}$ and then
\begin{equation}\label{3.6}
\int_{\Omega} [v_n(-\Delta)^\alpha\xi+h(v_n)\xi]dx=\int_{\Omega} \xi d\nu_{n}.
\end{equation}
By Lemma \ref{lm 2.1},
$$|(-\Delta)^\alpha\xi(x)|\le \frac{c_{9}\norm{\xi}_{L^\infty(\Omega)}}{1+|x|^{N+2\alpha}},\quad x\in \Omega.$$
Thus,
\begin{equation}\label{3.7}
 \lim_{n\to\infty}\int_{\Omega} v_n(x)(-\Delta)^\alpha\xi(x) dx=\int_{\Omega} u_\nu(x)(-\Delta)^\alpha\xi(x) dx.
\end{equation}

By (\ref{2.1.1}) and increasing monotonicity of $h$, it follows
$h(u_\nu)\le h(\mathbb{G}_{\R^N}[\nu])$ a.e. in $\R^N$.
By Lemma \ref{lm 2.2}, we find that $h(u_\nu)\in L^1(\R^N,d\omega)$. As a result,
for any $n\ge N_0$,
$$h(v_n)\to h(u_\nu)\quad{\rm in}\quad L^1(\R^N,d\omega).$$
Consequently,
\begin{equation}\label{3.8}
 \lim_{n\to\infty}\int_{\Omega} h(v_n)\xi(x) dx=\int_{\Omega} h(u_\nu)\xi(x) dx.
\end{equation}
It is obvious that
\begin{equation}\label{3.9}
 \lim_{n\to\infty}\int_{\Omega} \xi(x)d\nu_{n}=\int_{\Omega} \xi(x) d\nu.
\end{equation}
Combining (\ref{3.7}), (\ref{3.8}) with (\ref{3.9}) and taking $n\to\infty$ in (\ref{3.6}), we obtain that
\begin{equation}\label{3.10}
\int_{\Omega} [u_\nu(-\Delta)^\alpha\xi+h(u_\nu)\xi]dx=\int_{\Omega} \xi d\nu.
\end{equation}
Since $\xi\in \mathbb{X}_{\Omega}$ is arbitrary,  $u_\nu$
is a weak solution of (\ref{eq1.1}).\qquad$\Box$

\setcounter{equation}{0}
\section{Properties of weak solutions}

In this section, we investigate problem (\ref{eq1.3}). First we show that there is a unique solution of problem (\ref{eq1.3}), then we establish the
asymptotic behavior at the origin and infinity for the solution. In other words, we will prove
Theorem \ref{cor 1} and Theorem \ref{th 1}.

\subsection {Properties of $u_k$}

In this subsection, we consider the properties of nonnegative weak solution to (\ref{eq1.3}).
To this end, we introduce an auxiliary lemma.
\begin{lemma}\label{lm 2----1}\cite[Lemma 3.1]{CV3}
Assume that $v\in C^{2\alpha+\epsilon}(\bar B_1)$ with $\epsilon>0$
 satisfies
 $$(-\Delta)^\alpha v=\varphi\quad {\rm in}\quad B_1(0),$$
 where $\varphi\in C^1(\bar B_1)$. Then for $\beta\in (0,2\alpha)$, there exists $c_{11}>0$ such that
$$
\|v\|_{C^\beta(\bar B_{1/4}(0))}\le {c_{11}}(\|v\|_{L^\infty(B_1(0))}+\|\varphi\|_{L^\infty(B_1(0))}+\|v\|_{L^1(\R^N,d\omega)}).
$$

\end{lemma}

\bigskip

\begin{lemma}\label{lm 3.1}
Assume that $k>0$, $p\in(0,\frac{N}{N-2\alpha})$ and $u$ is a nonnegative weak solution of (\ref{eq1.3}).
Then  $u$ is a classical solution of (\ref{eq1.01}) and  for any $R>1$, there exists a weak solution $u_{R}$ of
\begin{equation}\label{3.11}
\arraycolsep=1pt
\begin{array}{lll}
 (-\Delta)^\alpha  u+u^p=k\delta_0\quad  & {\rm in}\quad B_R(0),\\[2mm]
  \phantom{(x----}
u=0 & {\rm in}\quad B_R^c(0)
\end{array}
\end{equation}
such that
\begin{equation}\label{3.1.7}
u-m_R\le u_{R}\le u\quad  \ {\rm in}\quad B_R(0)\setminus\{0\},
\end{equation}
where $m_R=\sup_{|x|>R}u(x)$.
\end{lemma}
{\bf Proof.} The proof is divide into two parts. First we show the regularity of solution $u$, then we find $u_R$ to
establish the inequality (\ref{3.1.7}).

{\it 1.Regularity of $u$.} Let $\{\eta_n\}\subset C^\infty_0(\R^N)$ be a sequence of radially decreasing and symmetric mollifiers such that supp$(\eta_n)\subset B_{\varepsilon_n}(0)$ with $\varepsilon_n\leq \frac1n$ and $u_n=u\ast\eta_n$.
We observe that
\begin{equation}\label{4.1}
u_n\to u\ \ {\rm and}\ \ u_n^p\to u^p\ \ {\rm in}\  \ L^1(\R^N,d\omega)\ \ {\rm as}\ n\to\infty.
\end{equation}
By Fourier transformation, we have that
$$\eta_n\ast(-\Delta)^{\alpha}\xi= (-\Delta)^{\alpha}(\xi\ast\eta_n),
$$
then
$$\int_{\R^N}\!\left(u(-\Delta)^{\alpha}(\xi\ast\eta_n) \!+\xi\ast\eta_nu^p\right)\!dx=\!
\int_{\R^N}\!\left(u\ast\eta_n(-\Delta)^{\alpha}\xi+\eta_n\ast u^p\xi\right)\! dx,
$$
where $\eta_n$ is radially  symmetric. It follows that $u_n$ is a classical solution of
\begin{equation}\label{eq1.4}
 \arraycolsep=1pt
\begin{array}{lll}
(-\Delta)^{\alpha}u_n+ u^p\ast\eta_n=k\eta_n\quad&\mbox{in }\ \ \R^N,\\[3mm]
\phantom{--}\lim_{|x|\to\infty}u_n(x)=0.
\end{array}
\end{equation}
We observe that $0\le u_n\le k\mathbb{G}_{\R^N}[\eta_n]$,
which implies $0\le u\le k\mathbb{G}_{\R^N}[\delta_0]$ in $\R^N\setminus\{0\}$.
Since $u^p\in L^1(\R^N,d\omega)$, we have $u^p\ast\eta_n\to u^p$ in $L^1(\R^N,d\omega)$ as $n\to\infty$
and that $\{k\eta_n+u_n^p-u^p\ast\eta_n\}$ converges to $k \delta_0$ in the distribution sense  as $n\to\infty$.

By Lemma \ref{teo CP}, we have that $0\le u_n\le  \mathbb{G}_{\R^N}[k\eta_n]$ and
 $\mathbb{G}_{\R^N}[k\eta_n]$ converges to $\mathbb{G}_{\R^N}[k\delta_0]$ uniformly in any compact set of $\R^N\setminus \{0\}$ and in $L^1(\R^N,d\omega)$. For a fixed $r>0$, there exists $N_0>0$ such that supp$(\eta_n)\subset \bar B_r(0)$ and there exists
 $c_{12}>0$ such that for any $n\ge N_0$,
$$\norm{u_n}_{L^\infty(B^c_{r/2}(0))}\le k\norm{\mathbb{G}_{\R^N}[\eta_n]}_{L^\infty(B^c_{r/2}(0))}\le c_{12}k\norm{\mathbb{G}_{\R^N}[\delta_0]}_{L^\infty( B_{r/2}^c(0))}$$  and $$ \norm{u_n}_{L^1(\R^N,d\omega)}\le k\norm{\mathbb{G}_{\R^N}[\eta_n]}_{L^1(\R^N,d\omega)}\le c_{12}k\norm{\mathbb{G}_{\R^N}[\delta_0]}_{L^1(\R^N,d\omega)}.$$
By Lemma \ref{lm 2----1}, for any $x_0\in \R^N$ with $|x_0|>4r$, there exists $\beta\in(0,2\alpha)$ such that
$$\arraycolsep=1pt
\begin{array}{lll}
\norm{u_n}_{C^{\beta}(B_{2r}(x_0))} \le {c_{11}}(\norm{u_n}_{L^1(\R^N,d\omega)}+\norm{u^p*\eta_n}_{L^{\infty}(B_{3r}(x_0))}+\norm{u_n}_{L^{\infty}(B_{3r}(x_0))})\\[2mm]
\phantom{-------} \le c_{11}(c_{12}k\|\mathbb{G}_{\R^N}[\delta_0]\|_{L^1(\R^N,d\omega)}+\norm{u^p}_{L^{\infty}(B_{3r}(x_0))}
 \\[2mm]\phantom{---------}+c_{12} k\norm{\mathbb{G}_{\R^N}[\delta_0]}_{L^\infty( B_{r/2}^c(0))})\\[2mm]
\phantom{-------} \le c_{11}(c_{12}k\|\mathbb{G}_{\R^N}[\delta_0]\|_{L^1(\R^N,d\omega)}+k^p\|\mathbb{G}_{\R^N}[\delta_0]\|^p_{L^{\infty}(B_{r/2}^c(0))}
 \\[2mm]\phantom{---------}+c_{12} k\norm{\mathbb{G}_{\R^N}[\delta_0]}_{L^\infty( B_{r/2}^c(0))}).
\end{array}
$$
Therefore,  by the definition of $u_n$ and  Arzela-Ascoli Theorem, we obtan that
$u\in C^{\frac\beta2}(B_{2r}(x_0))$. By Corollary 2.4 in \cite{RS}, we deduce that
$$\arraycolsep=1pt
\begin{array}{lll}
\norm{u_n}_{C^{2\alpha+\frac\beta2}(B_r(x_0))} \le {c_{13}}(\norm{u_n}_{L^1(\R^N,d\omega)}+\norm{u^p*\eta_n}_{C^{\frac\beta2}(B_{2r}(x_0))}
\\[2mm]\phantom{----------}+\norm{u_n}_{C^{\frac\beta2}(B_{2r}(x_0))})\\[2mm]
\phantom{--------} \le c_{14}(k\|\mathbb{G}_{\R^N}[\delta_0]\|_{L^1(\R^N,d\omega)}+\norm{u}_{C^{\frac\beta2}(B_{2r}(x_0))}
 \\[2mm]\phantom{----------}+
k^p\|\mathbb{G}_{\R^N}[\delta_0]\|^p_{L^{\infty}(B_{r/2}^c(0))}
+k\norm{\mathbb{G}_{\R^N}[\delta_0]}_{L^\infty( B_{r/2}^c(0))}),
\end{array}
$$
where $c_{13},c_{14}>0$. Thus, $u\in C^{2\alpha+\frac\beta4}(B_r(x_0))$ and by arbitrary of $r>0$ and $x_0$, $u$ is $C^{2\alpha+\frac\beta4}$ locally in $\R^N\setminus\{0\}$.
Therefore, $u_n\to u$ and $\eta_n\to 0$ uniformly in any compact subset of $\R^N\setminus\{0\}$ as $n\to\infty$.
We conclude  that $u$ is a classical solution of (\ref{eq1.01}) by Corollary 4.6 in \cite{CS1} .

{\it 2. Existence of $u_{R}$.}
It infers from (\ref{eq1.4}) that for given $R>1$,
 $u_n$ is a classical solution of
$$\arraycolsep=1pt
\begin{array}{lll}
(-\Delta)^{\alpha}u_n+ u_n^p=k\eta_n+u_n^p-u^p\ast\eta_n\qquad&\mbox{in }\ B_R(0),\\
\phantom{------}u_n\ge 0\quad &\mbox{in }\ B_R^c(0).
\end{array}
$$
We observe that for $R>1$,
$$\tilde u_n:=(u-m_{R-\epsilon_n})\ast\eta_n=u_n-m_{R-\epsilon_n}\le u_n\quad{\rm in}\quad \R^N$$ and
$(-\Delta)^{\alpha}\tilde u_n=(-\Delta)^{\alpha}u_n$, therefore,
$$(-\Delta)^{\alpha}\tilde u_n+ |\tilde u_n|^{p-1}\tilde u_n=k\eta_n+|\tilde u_n|^{p-1}\tilde u_n-u^p\ast\eta_n\quad{\rm in}\quad  B_R(0).$$
By the definition of $m_{R-\epsilon_n}$, we have $u-m_{R-\epsilon_n}\le 0$  in $B_{R-\epsilon_n}^c(0)$, and then
$$\tilde u_n\le 0\quad{\rm in}\  B_R^c(0).$$

Let $u_{n,R}$ be the solution of
$$\arraycolsep=1pt
\begin{array}{lll}
(-\Delta)^{\alpha}u_{n,R}+ u_{n,R}^p=k\eta_n+u_n^p-u^p\ast \eta_n\qquad&\mbox{in }\ B_R(0),\\
\phantom{------\  }u_{n,R}=0\quad &\mbox{in }\ B_R^c(0).
\end{array}
$$
By Lemma \ref{teo CP}, we have that
\begin{equation}\label{4.2}
\tilde u_n\le u_{n,R}\le u_n \quad {\rm in}\quad \R^N.
\end{equation}

It is known that $u_R:=\lim_{n\to\infty} u_{n,R}$ is a weak solution of (\ref{3.11}), since $\{k\eta_n+u_n^p-u^p\ast\eta_n\}$ converges to $k \delta_0$ in the distribution sense  as $n\to\infty$. Hence, (\ref{3.1.7}) follows by (\ref{4.2}) and $\tilde u_n\to u-m_R$  in $\R^N$ as $n\to\infty$. \qquad$\Box$

\medskip
With the help of Lemma \ref{lm 3.1}, we show next the uniqueness of weak solution to (\ref{eq1.3}).
\begin{proposition}\label{pr 3.1}
Assume that  $k>0$ and  $0<p<\frac{N}{N-2\alpha}$. Then (\ref{eq1.3}) admits a unique weak solution $u_k$.
\end{proposition}
{\bf Proof.}
{\it Existence.}
By Theorem \ref{teo 1}, there exists at least one weak solution $u_k$ to
$$(-\Delta)^\alpha  u+u^p=k\delta_0\quad {\rm in}\ \R^N$$
 such that $0\le u_k\le k\mathbb{G}_{\R^N}[\delta_0]$ a.e. in $\R^N$.
We observe that
 $$\mathbb{G}_{\R^N}[\delta_0](x)=\frac{c_{N,\alpha}}{|x|^{N-2\alpha}},\quad x\in\R^N\setminus\{0\},$$
then $$\lim_{|x|\to\infty}u_k(x)=0.$$
Thus $u_k$ is a weak solution of (\ref{eq1.3}).

{\it Uniqueness.}
We assume that  $u_k$, $\tilde u_k$ are two different weak solutions of (\ref{eq1.3}) and
$$A_0:=\min\{1,\limsup_{x\to0}|\tilde u_k-u_k|(x)\}.$$
We claim that $A_0>0$. In fact, if not, then $\lim_{x\to0}|\tilde u_k-u_k|(x)=0$.
Now we may assume that there exists $x_0\in\R^N\setminus\{0\}$ such that
$$(\tilde u_k-u_k)(x_0)=\sup_{x\in\R^N\setminus\{0\}}(\tilde u_k-u_k)(x)>0,$$
which implies that
$$(-\Delta)^\alpha(\tilde u_k-u_k)(x_0)\ge0.$$
Then we  obtain a contradiction by the fact that $\tilde u_k$ and $u_k$ are classical solutions of (\ref{eq1.01}) by Lemma \ref{lm 3.1}.
Therefore, $A_0>0$. Since
$$\lim_{|x|\to\infty}u_k(x)=0\quad{\rm and} \quad\lim_{|x|\to\infty}\tilde u_k(x)=0,$$
for $R>0$ large enough,
$$\epsilon_1:=\sup_{|x|\ge R}u_k(x)\le \frac{A_0}2\quad {\rm and}\quad\epsilon_2:=\sup_{|x|\ge R}\tilde u_k(x)\le \frac{A_0}2.$$
Since $u_k$ and $\tilde u_k$ are weak solutions of (\ref{eq1.3}), by Lemma \ref{lm 3.1},
there exist weak solutions $u_{k,R}$ and $\tilde u_{k,R}$ to (\ref{3.11}) such that
 \begin{equation}\label{3.1.2}
u_k-\epsilon_1\le  u_{k,R}\le u_k\quad {\rm in}\quad B_R(0)\setminus\{0\}
 \end{equation}
 and
 \begin{equation}\label{3.1.1}
\tilde u_k-\epsilon_2\le \tilde u_{k,R}\le \tilde u_k\quad {\rm in}\quad B_R(0)\setminus\{0\}.
 \end{equation}
Moreover, by Proposition \ref{pr 2.1} we obtain
$$  u_{k,R}\equiv\tilde  u_{k,R},$$
which, together with (\ref{3.1.1}) and (\ref{3.1.2}), implies that
$$|u_k-\tilde u_k|\le \max\{\epsilon_1,\epsilon_2\} \quad {\rm in }\quad B_R(0)\setminus\{0\}.$$
Thus,
$$\limsup_{x\to0}|u_k-\tilde u_k|(x)\le \max\{\epsilon_1,\epsilon_2\}<A_0.$$
This contradicts to the definition of $A_0$.
As a consequence,  problem (\ref{eq1.3}) has a unique weak solution. \qquad$\Box$

\medskip
Now we estimate the singularity rate of weak solution to (\ref{eq1.3}) at the origin.
\begin{proposition}\label{pr 3.2}
Let $k>0$,  $0<p<\frac{N}{N-2\alpha}$ and $u_k$  be the weak solution of (\ref{eq1.3}).
Then
\begin{equation}\label{3.1.3}
  \lim_{x\to0}u_k(x)|x|^{N-2\alpha}=c_{N,\alpha}k.
\end{equation}

\end{proposition}
{\bf Proof.}  On the one hand, we have that
\begin{equation}\label{eq 3.1--1}
  u_k(x)\le \mathbb{G}_{\R^N}[k\delta_0](x)=c_{N,\alpha}k|x|^{-N+2\alpha},\quad x\in\R^N\setminus\{0\}.
\end{equation}
On the other hand, from the proof of Theorem \ref{teo 1} and uniqueness of weak solution to (\ref{eq1.3}), we know
 $u_k=\lim_{R\to\infty} u_{k,R},$
 where $u_{k,R}$ is the weak solution of (\ref{3.11}).
By \cite[Lemma 2.1]{CV3} and \cite[Proposition1.1]{CV3},  we have that
$$\lim_{x\to0}u_{k,R}(x)|x|^{N-2\alpha}=c_{N,\alpha}k.$$
Then together with (\ref{eq 3.1--1}) and the fact that $\{u_{k,R}\}_R$ is an increasing sequence of functions,
(\ref{3.1.3}) holds.\qquad$\Box$
\medskip
\bigskip

\noindent{\bf Proof of Theorem \ref{cor 1}.}
By Proposition \ref{pr 3.1}, Lemma \ref{lm 3.1} and (\ref{3.1.3}), the assertion of Theorem \ref{cor 1} holds except part $(ii)$.

Now, we prove part $(ii)$ of Theorem \ref{cor 1}.
In fact, let $k_1\le k_2$ and $u_{k_1}$, $u_{k_2}$ be the solution of
(\ref{eq1.3}) with $k=k_1$ and $k=k_2$, respectively.
For $R>1$, we denote by $u_{k_1,R}$ and $u_{k_2,R}$ the solutions of
(\ref{3.11}) with $k=k_1$ and $k=k_2$, respectively. By $k_1\le k_2$ and Proposition \ref{pr 2.1}, we have that
$$u_{k_1,R}\le u_{k_2,R}\quad {\rm in }\ \ \R^N\setminus\{0\}.$$
Similar to the proof of Theorem \ref{teo 1}, we know that $u_{k_i}=\lim_{R\to\infty} u_{k_i,R}$ with $i=1,2$.
Therefore, $u_{k_1}\le u_{k_2}$ in $\R^N\setminus\{0\}$.
 \hfill$\Box$

\subsection{Asymptotic behavior of $u_1$ at $\infty$}

This subsection is devoted to investigate the asymptotic behavior of weak solution $u_1$ at $\infty$ to
\begin{equation}\label{eq 5.1}
 \arraycolsep=1pt
\begin{array}{lll}
 (-\Delta)^\alpha  u+u^p=\delta_0\quad   {\rm in}\quad \R^N,\\[2mm]
  \phantom{    }
  \lim_{|x|\to+\infty}u(x)=0,
\end{array}
\end{equation}
where  $p\in(1,\frac{N}{N-2\alpha})$.
We observe that
$$\lim_{x\to0}u_1(x)|x|^{N-2\alpha}=c_{N,\alpha}$$
and $u_1$ is a classical solution of
\begin{equation}\label{eq 5.2}
 (-\Delta)^\alpha  u+u^p=0\quad  {\rm in}\quad\R^N\setminus\{0\}.
\end{equation}


 In order to prove Theorem \ref{th 1}, we introduce some auxiliary lemmas.
For $\tau\in(-\infty,-N+2\alpha)$, we denote by $w_\tau$  a $C^2$ nonnegative radially symmetric function in $\R^N$ such that $w_\tau$ is decreasing in $|x|$ and for $|x|>1$,
\begin{equation}\label{6.4}
w_\tau(x)=\left\{
\arraycolsep=1pt
\begin{array}{lll}
|x|^\tau   &{\rm for}\quad \tau\in(-\infty,-N+2\alpha)\setminus\{-N\},\\[2mm]
  \phantom{    }
|x|^\tau\log^{\gamma_0}(\varrho_0+|x|) \quad&{\rm for}\quad \tau=-N,
\end{array}
 \right.
 \end{equation}
where $\varrho_0=e^{\frac1{2\alpha}}$ and $\gamma_0=\frac{N}{2\alpha}$.

\begin{lemma}\label{lm 5.1}
Assume that  $\tau\in(-\infty,-N+2\alpha)$.Then \\
$(i)$ for $\tau\in(-\infty,-N)$, there exist $R\ge 4$ and $c_{15}>1$ such that for $|x|>R$,
\begin{equation}\label{5.3}
  \frac1{c_{15}}|x|^{-N-2\alpha}\le -(-\Delta)^\alpha w_{\tau}(x)\le c_{15}|x|^{-N-2\alpha};
\end{equation}
$(ii)$ for $\tau=-N$, there exist $R\ge 4$ and $c_{15}>1$ such that for $|x|>R$,
\begin{equation}\label{5.4}
 \frac1{c_{15}}|x|^{-N-2\alpha}\log^{\gamma_0+1}|x|\le- (-\Delta)^\alpha w_{\tau}(x)\le c_{15}|x|^{-N-2\alpha}\log^{\gamma_0+1}|x|;
\end{equation}
$(iii)$ for $\tau\in(-N,-N+2\alpha)$, there exist $R\ge4$ and $c_{15}>1$ such that for $|x|>R$,
\begin{equation}\label{5.05}
 \frac1{c_{15}}|x|^{\tau-2\alpha}\le- (-\Delta)^\alpha w_{\tau}(x)\le c_{15}|x|^{\tau-2\alpha}.
\end{equation}

\end{lemma}
{\bf Proof.}  In the following, we shall use the equivalent definition of $(-\Delta)^\alpha w_\tau$, that is,
$$(-\Delta)^\alpha w_\tau(x)=-\frac12\int_{\R^N}\frac{w_\tau(x+z)+w_\tau(x-z)-2w_\tau(x)}{|z|^{N+2\alpha}}dz.$$
$(i)$ The case of $\tau\in(-\infty,-N)$.
On the one hand,  for $|x|>4$, we have that
\begin{equation}\label{5.5}
\arraycolsep=1pt
\begin{array}{lll}
 -(-\Delta)^\alpha w_{\tau }(x)=\frac12\int_{\R^N\setminus (B_1(x)\cup B_1(-x))}\frac{w_\tau(x+z)+w_\tau(x-z)-2w_\tau(x)}{|z|^{N+2\alpha}}dz
 \\[3mm]\phantom{-------\ }
 +\frac12\int_{B_1(x)\cup B_1(-x)}\frac{w_\tau(x+z)+w_\tau(x-z)-2w_\tau(x)}{|z|^{N+2\alpha}}dz
 \\[3mm]\phantom{------\ \ }
 \le\frac{|x|^{\tau -2\alpha}}2\int_{D_0}\frac{I_x(y)}{|y|^{N+2\alpha}}dy+c_{16}|x|^{-N-2\alpha}
\end{array}
 \end{equation}
where $c_{16}>0$ depends on $\norm{w}_{L^1(B_1(0))}$, $e_x=\frac{x}{|x|}$, $D_0=\R^N\setminus (B_{\frac1{|x|}}(e_x)\cup B_{\frac1{|x|}}(-e_x))$ and
$$I_x(y) = |e_x+y|^{\tau }+|e_x-y|^{\tau }-2.$$
On the other hand, for $|x|\ge 4$,
\begin{equation}\label{5.6}
\arraycolsep=1pt
\begin{array}{lll}
 -(-\Delta)^\alpha w_{\tau }(x)=\frac12\int_{\R^N\setminus (B_1(x)\cup B_1(-x))}\frac{w_\tau(x+z)+w_\tau(x-z)-2w_\tau(x)}{|z|^{N+2\alpha}}dz
 \\[3mm]\phantom{-------\ }
 +\frac12\int_{B_1(x)\cup B_1(-x)}\frac{w_\tau(x+z)+w_\tau(x-z)-2w_\tau(x)}{|z|^{N+2\alpha}}dz
 \\[3mm]\phantom{------\ \ } \ge\frac{|x|^{\tau -2\alpha}}2\int_{D_0}\frac{I_x(y)}{|y|^{N+2\alpha}}dy-2\int_{B_1(x)}\frac{w_\tau(x)}{|z|^{N+2\alpha}}dz
 \\[3mm]\phantom{------\ \ }\ge \frac{|x|^{\tau -2\alpha}}2\int_{D_0}\frac{I_x(y)}{|y|^{N+2\alpha}}dy-c_{17}|x|^{\tau-N-2\alpha},
\end{array}
 \end{equation}
where $c_{17}>0$.

 \smallskip
\emph{Claim 1.  There exists $c_{18}>1$ such that
 \begin{equation}\label{15-09-1}
 \arraycolsep=1pt
 \begin{array}{lll}
  \frac1{c_{18}}|x|^{-N-\tau }\le \int_{D_1\cup D_2}\frac{I_x(y)}{|y|^{N+2\alpha}}dy\le c_{18} |x|^{-N-\tau },
 \end{array}
 \end{equation}
where $D_1=B_{\frac12}(-e_x)\setminus B_{\frac1{|x|}}(-e_x)$ and $D_2=B_{\frac12}(e_x)\setminus B_{\frac1{|x|}}(e_x)$.}

In fact, for $y\in D_1$, we observe that
$$-2\le |e_x-y|^{\tau }-2 \le -1\quad {\rm and }\quad \frac12\le |y|\le \frac32,$$
 then
\begin{eqnarray*}
  \int_{D_1}\frac{I_x(y)}{|y|^{N+2\alpha}}dy &\le& c_{19}\int_{B_{\frac12}(0)\setminus B_{\frac1{|x|}}(0)}|y|^{\tau }dy+\int_{D_1}\frac{|e_x-y|^{\tau }-2}{|y|^{N+2\alpha}}dy
 \\&\le& c_{20}\int_{|x|^{-1}}^{\frac12}r^{\tau +N-1}dr
   \\& \le& c_{21} |x|^{-N-\tau }
\end{eqnarray*}
and
\begin{eqnarray*}
  \int_{D_1}\frac{I_x(y)}{|y|^{N+2\alpha}}dy &\ge&
  c_{22}\int_{|x|^{-1}}^{\frac12}r^{\tau +N-1}dr-c_{23}
   \\& \ge& c_{24}|x|^{-N-\tau}-c_{23},
\end{eqnarray*}
where $c_{18},...,c_{24}$ are positive constants. Since $-N-\tau>0$, there exist $R\ge4$ and $c_{25}>0$ such that
for $|x|\ge R$,
 $$ \frac1{c_{25}}|x|^{-N-\tau }\le \int_{D_1}\frac{I_x(y)}{|y|^{N+2\alpha}}dy\le c_{25} |x|^{-N-\tau }.$$
By the fact
$$\int_{D_1}\frac{I_x(y)}{|y|^{N+2\alpha}}dy=\int_{D_2}\frac{I_x(y)}{|y|^{N+2\alpha}}dy,$$
we obtain (\ref{15-09-1}).
\smallskip

\emph{Claim 2.  There exists $c_{26}>0$ such that
 \begin{equation}\label{15-09-2}
  | \int_{B_{\frac12}(0)}\frac{I_x(y)}{|y|^{N+2\alpha}}dy|\le c_{26}.
 \end{equation}}
Indeed, since function $I_x$ is $C^2$ in $\bar B_{\frac12}(0)$ such that
$$I_x(0)=0\quad{\rm and}\quad  I_x(y)=I_x(-y),\quad y\in \bar B_{\frac12}(0),$$
then $\nabla I_x(0)=0$ and there exists $c_{27}>0$ such that
$$|D^2 I_x(y)|\le c_{27}, \quad y\in B_{\frac12}(0). $$
Then we have
$$I_x(y)\le c_{27}|y|^2,\quad y\in B_{\frac12}(0),$$
which implies that
$$ |\int_{B_{\frac12}(0)}\frac{I_x(y)}{|y|^{N+2\alpha}}dy|\le c_{27}\int_{B_{\frac12}(0)}\frac{|y|^{2}}{|y|^{N+2\alpha}}dy\le c_{26}.
$$

\emph{Claim 3.  There exists $c_{28}>0$ such that
 \begin{equation}\label{15-09-3}
   |\int_{D_3}\frac{I_x(y)}{|y|^{N+2\alpha}}dy|\le c_{28},
 \end{equation}
where $D_3=\R^N\setminus (B_{\frac12}(0)\cup B_{\frac12}(e_x)\cup B_{\frac12}(-e_x))=D_0\setminus(D_1\cup D_2\cup B_{\frac12}(0))$.}

In fact, for $y\in D_3$, we observe that there exists $c_{29}>0$ such that
$|I_x(y)|\le c_{29}$ and
$$
   |\int_{D_3}\frac{I_x(y)}{|y|^{N+2\alpha}}dy|\le \int_{\R^N\setminus B_{\frac12}(0)}\frac{c_{29}}{|y|^{N+2\alpha}}dy\le c_{30},
$$
where $c_{30}>0$.
Since $\lim_{|x|\to\infty}|x|^{-N-\tau }=\infty$ for $\tau<-N$, by (\ref{15-09-1})-(\ref{15-09-3}), there exist $R\ge 4$ and $c_{31}>1$ such that for $|x|\ge R$,
\begin{equation}
\arraycolsep=1pt
\begin{array}{lll}
\frac1{c_{31}}|x|^{-N-2\alpha}\le  -(-\Delta)^\alpha w(x)
\le c_{31} |x|^{-N-2\alpha}.
\end{array}
\end{equation}

$(ii)$ The case of $\tau=-N$.  Similarly to (\ref{5.5}) and (\ref{5.6}),
we have that for $|x|>4$,
$$
\arraycolsep=1pt
\begin{array}{lll}
 -(-\Delta)^\alpha w_{\tau }(x) \le\frac{|x|^{-N-2\alpha}\log^{\gamma_0}(\varrho_0+|x|)}2\int_{D_0}\frac{II_x(y)}{|y|^{N+2\alpha}}dy+c_{32}|x|^{-N-2\alpha}
\end{array}
$$
and
$$
\arraycolsep=1pt
\begin{array}{lll}
 -(-\Delta)^\alpha w_{\tau }(x) \ge\frac{|x|^{-N-2\alpha}\log^{\gamma_0}(\varrho+|x|)}2\int_{D_0}\frac{II_x(y)}{|y|^{N+2\alpha}}dy-
 c_{32}|x|^{-2N-2\alpha}\log^{\gamma_0}|x|,
\end{array}
$$
where $c_{32}>0$ and
$$\arraycolsep=1pt
\begin{array}{lll}
II_x(y) =\frac{\log^{\gamma_0}(\varrho_0+|x||e_x+y|)}{\log^{\gamma_0}(\varrho_0+|x|)}|e_x+y|^{-N }+\frac{\log^{\gamma_0}(\varrho_0+|x||e_x-y|)}{\log^{\gamma_0}(\varrho_0+|x|)}
|e_x-y|^{-N}-2
\end{array}.
$$
For $y\in D_1$, we have that
$\frac{\log^{\gamma_0}(\varrho_0+|x||e_x+y|)}{\log^{\gamma_0}(\varrho_0+|x|)}\le 1$,
 then
\begin{eqnarray*}
  \int_{D_1}\frac{II_x(y)}{|y|^{N+2\alpha}}dy &\le& c_{33}\int_{B_{\frac12}(0)\setminus B_{\frac1{|x|}}(0)}|y|^{-N}dy+c_{34}
   \\& \le& c_{35} \log|x|+c_{34}
\end{eqnarray*}
and
\begin{eqnarray*}
  \int_{D_1}\frac{II_x(y)}{|y|^{N+2\alpha}}dy &\ge&
  c_{33}\int_{|x|^{-1}}^{\frac12}r^{-1}\frac{\log^{\gamma_0}(\varrho_0+|x|r)}{\log^{\gamma_0}(\varrho_0+|x|)}dr-c_{34}
  \\& =&   c_{33}\int_{1}^{\frac12|x|}s^{-1}\frac{\log^{\gamma_0}(\varrho_0+s)}{\log^{\gamma_0}(\varrho_0+|x|)}ds-c_{34}
   \\& \ge& c_{35}\log|x|-c_{34},
\end{eqnarray*}
where $c_{33},c_{34},c_{35}>0$.
By the fact that
$$\int_{D_1}\frac{II_x(y)}{|y|^{N+2\alpha}}dy=\int_{D_2}\frac{II_x(y)}{|y|^{N+2\alpha}}dy,$$
there exists $c_{36}>0$ such that
\begin{equation}\label{7.1}
 \arraycolsep=1pt
 \begin{array}{lll}
  \frac1{c_{36}}\log|x|\le \int_{D_1\cup D_2}\frac{II_x(y)}{|y|^{N+2\alpha}}dy\le c_{36} \log|x|.
 \end{array}
 \end{equation}
Similarly to  (\ref{15-09-1})-(\ref{15-09-3}), there exists $c_{37}>0$ such that
\begin{equation}\label{7.2}
 \arraycolsep=1pt
 \begin{array}{lll}
 |\int_{\R^N\setminus (B_{\frac12}(e_x)\cup B_{\frac12}(-e_x))}\frac{II_x(y)}{|y|^{N+2\alpha}}dy|\le c_{37},
 \end{array}
 \end{equation}
which, together with (\ref{7.1}), imply (\ref{5.4}).

\medskip
$(iii)$ The case that $\tau\in(-N,-N+2\alpha)$.
By Lemma 3.1 and Lemma 3.2 in \cite{FQ2}, we have
\begin{equation}\label{7.3}
  (-\Delta)^\alpha |x|^{\tau }=c(\tau)|x|^{\tau -2\alpha},
\end{equation}
where $c(\tau)<0$.

Let $\tilde w(x)=w_{\tau }(x)-|x|^{\tau }$ for $x\in \R^N\setminus\{0\}$,
then $\tilde w=0$ in $B_1^c(0)$.
For $|x|>4$, we have that
\begin{equation}\label{5.7}
\arraycolsep=1pt
\begin{array}{lll}
  |(-\Delta)^\alpha\tilde w(x)| \le\int_{B_1(0)}\frac{w_\tau(z)+|z|^{\tau }}{|z-x|^{N+2\alpha}}dz \\[3mm]
  \phantom{------}
   \le(|x|-1)^{-N-2\alpha}  \int_{B_1(0)}(w_\tau(z)+|z|^{\tau })dz,
\end{array}
\end{equation}
which, together with (\ref{7.3}), imply (\ref{5.05}). We complete the proof.\qquad$\Box$

\begin{lemma}\label{lm 5.2}
Let $\eta:\R^N\to[0,1]$ be a $C^2$ function with support in $B_2(0)$ and $\eta=1$ in $B_1(0)$,
$\bar w(x)=\eta(x)|x|^{-N+2\alpha}$ for $x\in \R^N$.
Then for $|x|>4$, there exists $c_{38}>1$ such that
\begin{equation}\label{5.8}
\frac1{c_{38}}|x|^{-N-2\alpha}\le -(-\Delta)^\alpha\bar w(x)\le c_{38}|x|^{-N-2\alpha},\quad x\in B_4^c(0).
\end{equation}
 \end{lemma}
{\bf Proof.} For $|x|>4$, we have that
$$
\arraycolsep=1pt
\begin{array}{lll}
  -(-\Delta)^\alpha\bar w(x) =\int_{B_2(0)}\frac{\bar w(z)}{|z-x|^{N+2\alpha}}dz  \le (|x|-2)^{-N-2\alpha}  \int_{B_2(0)}\bar w(z)dz
\end{array}
$$
and
$$
\arraycolsep=1pt
\begin{array}{lll}
  -(-\Delta)^\alpha\bar w(x)  \ge (|x|+2)^{-N-2\alpha}  \int_{B_2(0)}\bar w(z)dz,
\end{array}
$$
which, together with $\int_{B_2(0)}\bar w(z)dz<+\infty$, imply (\ref{5.8}). \qquad$\Box$

\medskip

Now we are in the position to prove Theorem \ref{th 1}.
\smallskip

\noindent{\bf Proof of Theorem \ref{th 1}.} For $p\in(1,\frac{N}{N-2\alpha})$,
we denote
\begin{equation}\label{tau_p}
\tau_p=\left\{
\arraycolsep=1pt
\begin{array}{lll}
 -\frac{2\alpha}{p-1}\qquad   &{\rm for}\quad p\in[1+\frac{2\alpha}{N},\frac{N}{N-2\alpha}),\\[2mm]
  \phantom{    }
 -\frac{N+2\alpha}{p}&{\rm for}\quad p\in(1,1+\frac{2\alpha}{N}).
\end{array}
 \right.
 \end{equation}
 We note that $\tau_p$ is continuous and strictly increasing with respect to $p$,  $\tau_p=-N$ if $p=1+\frac{2\alpha}{N}$ and
 $\lim_{p\to0^+}\tau_p=-\infty$.

{\it Lower bound.} Since $\lim_{|x|\to0}u_1(x)=\infty$ and $u_1$ is continuous and positive in $\R^N\setminus\{0\}$,
then there exists $c_{39}>0$ such that
\begin{equation}\label{5.9}
 c_{39}w_{\tau_p}\le u_1\quad {\rm in}\ \ \bar B_R(0)\setminus\{0\},
\end{equation}
where $R>4$ is from Lemma \ref{lm 5.1}.

We note that for $p\in(1,1+\frac{2\alpha}{N})$, $\tau_pp=-N-2\alpha$;
for $p=1+\frac{2\alpha}{N}$, $\tau_p=-N$, $\gamma_0+1=p\gamma_0$ and
for $p\in(1+\frac{2\alpha}{N},\frac{N}{N-2\alpha})$, $\tau_p-2\alpha=p\tau_p$. By Lemma \ref{lm 5.1}, there
exists $t_0\in(0,1)$ such that
$$(-\Delta)^\alpha (t_0c_{39}w_{\tau_p})+(t_0c_{39}w_{\tau_p})^p\le 0\quad {\rm in }\quad B_R^c(0).$$

We claim that $u_1\ge t_0c_{39}w_{\tau_p}$ in $B_R^c(0)$.
In fact, if not, there would exist $x_0\in B_R^c$ such that
\begin{eqnarray*}
 (u_1-t_0c_{39}w_{\tau_p})(x_0) &=& \min_{x\in B_R^c(0)}(u_1-t_0c_{39}w_{\tau_p})(x) \\
   &=&\min_{x\in \R^N\setminus\{0\}}(u_1-t_0c_{39}w_{\tau_p})(x)<0,
\end{eqnarray*}
since $u_1-t_0c_{39}w_{\tau_p}\ge0$ in $\bar B_R(0)$ and $\lim_{|x|\to\infty}(u_1-t_0c_{39}w_{\tau_p})(x)=0$.
Then $(-\Delta)^\alpha(u_1-t_0c_{39}w_{\tau_p})(x_0)<0$. However, $$
       (-\Delta)^\alpha(u_1-t_0c_{39}w_{\tau_p})(x_0) \ge -u_1^p(x_0)+(t_0c_{39}w_{\tau_p})^p(x_0)>0,$$
which is a contradiction.

\smallskip
{\it Upper bound.} Since $\lim_{x\to0}u_1(x)|x|^{N-2\alpha}=c_{N,\alpha}$,
there exists $c_{40}>0$ such that
$u_1(x)\le c_{40}|x|^{-N+2\alpha}$ in $B_1(0)\setminus\{0\}$. Then
there exists $c_{41}>1$ such that
$$u_1\le c_{40}w_{\tau_p}+c_{41}\bar w\quad{\rm in }\ \ \bar B_R(0),$$
where $\bar w$ is from Lemma \ref{lm 5.2}.
Denote by $W=c_{40}w_{\tau_p}+c_{41}\bar w$.  By Lemma \ref{lm 5.1} and Lemma \ref{lm 5.2}, there
exist $t_1>1$ such that
$$(-\Delta)^\alpha (t_1W)+(t_1W)^p\ge 0\quad {\rm in }\quad B_R^c(0).$$

We claim that $u_1\le t_1W$ in $B_R^c(0)$.
In fact, if not, then there exists $x_1\in B_R^c(0)$ such that
\begin{eqnarray*}
 (u_1-t_1W)(x_1) &=& \max_{x\in B_R^c(0)}(u_1-t_1W)(x) \\
   &=& \max_{x\in \R^N\setminus\{0\}}(u_1-t_1W)(x)>0.
\end{eqnarray*}
Thus, $(-\Delta)^\alpha(u_1-t_1W)(x_1)>0$. But
$$
       (-\Delta)^\alpha(u_1-t_1W)(x_1) \le -u_1^p(x_1)+(t_1W)^p(x_1)<0,$$
we obtain a contradiction.

Since $\bar w=0$ in $B_2^c$, combining (\ref{6.4}) with (\ref{tau_p}), we obtain the decays of $u_1$ for $p\in(1,\frac{N}{N-2\alpha})$.\qquad$\Box$

\setcounter{equation}{0}
\section {Properties of the limit function}

Let $u_\infty$ be given by (\ref{definition infty}) and $u_{k,\Omega}$ be a weak solution of (\ref{eq1.001})
when $\Omega$ is an unbounded regular domain including the origin. We plan to study properties of both $u_\infty$ and $u_{k,\Omega}$.

\subsection{Properties of $u_\infty$ }
This subsection is devoted to prove Theorem \ref{teo 2}. To this end, we introduce some propositions.

\begin{proposition}\label{lm 13-09-1}
Assume that $p\in(1,\frac{N}{N-2\alpha})$ and $u_\infty$ is defined  in (\ref{definition infty}).
Then
\begin{equation}\label{3.2.1}
u_\infty(x)=|x|^{-\frac{2\alpha}{p-1}}u_\infty(\frac{x}{|x|}),\quad x\in\R^N\setminus\{0\}.
\end{equation}
\end{proposition}
{\bf Proof.} For $\lambda>0$, we denote
$$\tilde u_\lambda(x)=\lambda^{\frac{2\alpha}{p-1}}u_k(\lambda x),\quad x\in \R^N\setminus\{0\}, $$
where $u_k$ is the solution of (\ref{eq1.3}).
By direct computation, we have for $x\in\R^N\setminus\{0\}$ that,
\begin{eqnarray}
 (-\Delta)^\alpha  \tilde u_\lambda(x) +\tilde u_\lambda^p(x)
  &=&\lambda^{\frac{2\alpha p}{p-1}}[(-\Delta)^\alpha u_k(\lambda x)  +  u_k^p(\lambda x)] \nonumber \\
  &=&0.\label{13-09-5}
\end{eqnarray}
Moreover, for $f\in C_0(\R^N)$,
\begin{eqnarray*}
 \langle(-\Delta)^\alpha  \tilde u_\lambda +\tilde u_\lambda^p, f\rangle&=&\lambda^{\frac{2\alpha p}{p-1}}\int_{\R^N}
 [(-\Delta)^\alpha u_k(\lambda x)  +  u_k^p(\lambda x)]f(x)dx\nonumber
  \\&=&\lambda^{\frac{2\alpha p}{p-1}-N}\int_{\R^N}  [(-\Delta)^\alpha u_k(z)  +  u_k^p(z)]f(\frac{z}{\lambda})dz\nonumber
   \\  &=&\lambda^{\frac{2\alpha p}{p-1}-N}kf(0),
\end{eqnarray*}
where $\frac{2\alpha p}{p-1}-N>0$ by the fact that $p\in(1,\frac{N}{N-2\alpha})$.
Thus,
\begin{equation}\label{3.2.2}
  (-\Delta)^\alpha  \tilde u_\lambda +\tilde u_\lambda^p=\lambda^{\frac{2\alpha p}{p-1}-N}k\delta_0\quad {\rm in}\ \ \R^N.
\end{equation}
We observe that $\lim_{|x|\to\infty}\tilde u_\lambda(x)=0$
and  $u_{k\lambda^{\frac{2\alpha p}{p-1}-N}}$ is a unique weak solution of
(\ref{eq1.3}) with $k$ replaced by $\lambda^{\frac{2\alpha p}{p-1}-N}k$,
then  for $x\in \R^N\setminus\{0\}$,
\begin{equation}\label{21-10-1}
u_{k\lambda^{\frac{2\alpha p}{p-1}-N}}(x)=\tilde u_\lambda (x)=\lambda^{\frac{2\alpha }{p-1}}u_k(\lambda x)
\end{equation}
and letting $k\to\infty$ we have that
$$
u_{\infty}(x)=\lambda^{\frac{2\alpha }{p-1}}u_\infty(\lambda x),\qquad x\in \R^N\setminus\{0\},
$$
which implies (\ref{3.2.1}) by taking $\lambda=|x|^{-1}$.
\qquad$\Box$ \medskip

\begin{proposition}\label{pr 3.3}
 Suppose that $p\in(0,1+\frac{2\alpha}{N}]$ and  $u_\infty$ is given by (\ref{definition infty}).
Then $$u_\infty=\infty\quad{\rm in}\quad \R^N.$$

\end{proposition}
{\bf Proof.}  {\it In the case of $p\in(0,1]$.}  We observe that $$\mathbb{G}_{\R^N}[\delta_0],\,\quad \mathbb{G}_{\R^N}[(\mathbb{G}_{\R^N}[\delta_0])^p]>0$$ in $\R^N$.
Since $$u_k\ge k\mathbb{G}_{\R^N}[\delta_0]-k^p\mathbb{G}_{\R^N}[(\mathbb{G}_{\R^N}[\delta_0])^p],$$
we obtain  $\lim_{k\to\infty}u_k=\infty$ in $\R^N$ for  $p\in(0,1)$.
For $p=1$, we see that $u_k=ku_1$. Hence, $\lim_{k\to\infty}u_k=\infty$ in $\R^N$ by the fact that $u_1>0$ in $\R^N$.

\smallskip
{\it In the case of $p\in(1,1+\frac{2\alpha}{N}]$.}
It derives from (\ref{3.2.1}) that
$$u_\infty(x)\ge c_{42}|x|^{-\frac{2\alpha}{p-1}},\quad x\in\R^N\setminus\{0\},$$
where $c_{42}=\min_{|x|=1}u_\infty(x)>0$, since $u_\infty\ge u_{k}$ in $\R^N\setminus\{0\}$.
Since $u_\infty=\lim_{k\to\infty} u_k$ in $\R^N\setminus\{0\}$,
we deduce that
\begin{equation}\label{3.2.3}
\pi_k:=\int_{B_{\frac14}(0)}u_k(x)dx\to\infty\quad {\rm as}\ k\to\infty.
\end{equation}

Fix $y_0\in \R^N$ such that $|y_0|=1$, it follows by Lemma \ref{lm 2---1} that problem
\begin{equation}\label{5.1}
\arraycolsep=1pt
\begin{array}{lll}
 (-\Delta)^\alpha  u+u^p=0 \quad & {\rm in}\quad  B_{\frac14}(y_0),\\[2mm]
 \phantom{  (-\Delta)^\alpha  +u^p}
u=0  \quad & {\rm in}\quad \R^N \setminus (B_{\frac14}(y_0)\cup B_{\frac14}(0)),\\[2mm]
\phantom{ (-\Delta)^\alpha  +u^p}
u=u_k  \quad & {\rm in}\quad B_{\frac14}(0)
\end{array}
\end{equation}
admits a unique solution $w_k$.
By Lemma \ref{teo CP},
\begin{equation}\label{4.1.3}
 u_{k}\ge w_k\quad {\rm in}\quad B_{\frac14}(y_0).
\end{equation}
Let  $\tilde w_k=w_k-u_k\chi_{B_{\frac14}(0)},$
then $\tilde w_k=w_k$ in $B_{\frac14}(y_0)$ and for $x\in B_{\frac14}(y_0)$,
$$
\arraycolsep=1pt
\begin{array}{lll}
(-\Delta)^\alpha \tilde w_k(x) =
-\lim_{\epsilon\to0^+}\int_{B_{\frac14}(y_0)\setminus B_\epsilon(x)}\frac{w_k(z)-w_k(x)}{|z-x|^{N+2\alpha}}dz
\\[3mm]\phantom{--------}+\lim_{\epsilon\to0^+}\int_{B_{\frac14}^c(y_0)\setminus B_\epsilon(x)}\frac{w_k(x)}{|z-x|^{N+2\alpha}}dz
\\[3mm]\phantom{------}
=-\lim_{\epsilon\to0^+}\int_{\R^N\setminus B_\epsilon(x)}\frac{w_k(z)-w_k(x)}{|z-x|^{N+2\alpha}}dz +\int_{B_{\frac14}(0)}\frac{u_k(z)}{|z-x|^{N+2\alpha}}dz
\\[3mm]\phantom{------}\ge(-\Delta)^\alpha w_k(x)+c_{43}\pi_k,
\end{array}
$$
where $c_{43}=(\frac45)^{N+2\alpha}$ and the last inequality follows by the fact of $$|z-x|\le |x|+|z|\le 5/4\quad
{\rm for}\ z\in B_{\frac14}(0), x\in B_{\frac14}(y_0).$$
Therefore,
\begin{eqnarray*}
(-\Delta)^\alpha \tilde w_k(x)+\tilde w_k^p(x) &\ge&  (-\Delta)^\alpha w_k(x)+w_k^p(x)+ c_{43}\pi_k \\
        &=&c_{43}\pi_k, \qquad x\in B_{\frac14}(y_0).
     \end{eqnarray*}
that is,  $\tilde w_k$  is a super solution of
\begin{equation}\label{4.1.2}
\arraycolsep=1pt
\begin{array}{lll}
\displaystyle (-\Delta)^\alpha  u+u^p=c_{43}\pi_k \quad & {\rm in}\quad  B_{\frac14}(y_0),\\[2mm]
 \phantom{  (-\Delta)^\alpha  +u^{p,}}
u=0  \quad & {\rm in}\quad B_{\frac14}^c(y_0).
\end{array}
\end{equation}
Let $\eta_1$ be the solution of
$$
\arraycolsep=1pt
\begin{array}{lll}
 (-\Delta)^\alpha  u=1 \quad & {\rm in}\quad  B_{\frac14}(y_0),\\[2mm]
 \phantom{  (-\Delta)^\alpha  }
u=0  \quad & {\rm in}\quad  B^c_{\frac14}(y_0).
\end{array}
$$
Then
 $(c_{43}\pi_k)^{\frac1p} \frac{\eta_1}{2\max_{\R^N}\eta_1}$ is sub solution of (\ref{4.1.2}) for $k $ large enough. By Lemma \ref{teo CP}, for $k$ big we have
 $$\tilde w_k(x)\ge (c_{43}\pi_k)^{\frac1p} \frac{\eta_1(x)}{2\max_{\R^N}\eta_1},\quad \forall  x\in B_{\frac14}(y_0),$$
 which implies that
 $$w_k(y_0)\ge \frac{(c_{43}\pi_k)^{\frac1p}}2.$$
 Therefore, (\ref{4.1.3}) and (\ref{3.2.3}) imply $$u_\infty(y_0)=\lim_{k\to\infty}u_{k}(y_0)\ge \lim_{k\to\infty}w_k(y_0)=\infty.$$
Since $y_0$ is arbitrary on $\partial B_1(0)$, by (\ref{3.2.1}), it follows that $u_\infty=\infty$ in $\R^N$.
$\Box$

\begin{proposition}\label{pr 3.4}
 Suppose that $p\in(1+\frac{2\alpha}N,\frac{N}{N-2\alpha})$ and  $u_\infty$ is given by (\ref{definition infty}).
Then $u_\infty$ is a classical solution of (\ref{eq1.01}).

\end{proposition}
{\bf Proof.} For  $p\in(1+\frac{2\alpha}N,\frac{N}{N-2\alpha})$, we observe that $\tau_p:=-\frac{2\alpha}{p-1}\in(-N,-N+2\alpha),$
$\tau_p-2\alpha=\tau_pp$
and $$(-\Delta)^\alpha |x|^{\tau_p}=c(\tau_p)|x|^{\tau_p-2\alpha},\quad x\in\R^N\setminus\{0\},$$
where $c(\tau_p)<0$, see Lemma 3.1 in \cite{FQ2}  and Lemma 3.2 in \cite{FQ2}. Let
$$W_p(x)=[-c(\tau_p)]^{\frac1{p-1}}|x|^{\tau_p}, \quad x\in\R^N\setminus\{0\}.$$
Then, $W_p$ is a solution of (\ref{eq1.01}).

{\it  We first prove that}
\begin{equation}\label{4.2.1}
 u_\infty\le W_p\quad{\rm in}\ \ \R^N\setminus\{0\}.
\end{equation}
 In fact, we observe that $u_k=\lim_{R\to\infty} u_{k,R}$ in $\R^N\setminus\{0\}$,
where $u_{k,R}$ is the solution of (\ref{3.11}) with $R>1$ and
$$\lim_{x\to0}u_{k,R}(x)|x|^{N-2\alpha}=c_{N,\alpha}k.$$
Then $$\lim_{x\to0}\frac{u_{k,R}(x)}{W_p(x)}=0.$$
Moreover, we know that $u_{k,R}$ is a classical solution of
$$
\arraycolsep=1pt
\begin{array}{lll}
(-\Delta)^\alpha u+u^p=0\quad {\rm in}\quad B_R(0)\setminus\{0\},
\\[2mm]\phantom{----\ \ }
\quad u=0\quad {\rm in}\quad B_R^c(0).
\end{array}
$$
By Lemma \ref{teo CP} with $O=B_R(0)\setminus B_\epsilon(0)$ and $\epsilon>0$ small enough, we obtain
that $$u_{k,R}\le W_p\quad{\rm in}\ \ \R^N\setminus\{0\},$$
which implies that for any $k>0$,
$$u_{k}\le W_p\quad{\rm in}\ \ \R^N\setminus\{0\}.$$
Thus, by the definition of $u_\infty$, (\ref{4.2.1}) holds.

{\it Next we prove that $u_\infty$ is a solution of (\ref{eq1.01}).} 
 We observe that $W_p\in L^1(\R^N, d\omega)$, then for  any $x_0\not=0$,
  there exist $c_{44},c_{45}>0$ independent of $k$ such that
$$\|u_k\|_{L^1(\R^N,d\omega)}\le c_{44}\quad{\rm and}\quad \|u_k\|_{L^\infty(B_{\frac{|x_0|}2}(x_0))}\le c_{45}.$$
It follows by the same argument in the proof of regularity in Lemma \ref{lm 3.1} that there  exist $\epsilon>0$ and ${c_{46}}>0$ independent of $k$ such that
$$\norm{u_k}_{C^{2\alpha+\epsilon}(B_{\frac{|x_0|}4}(x_0))}\le {c_{46}}.$$
By the definition of $u_\infty$ and the Arzela-Ascoli Theorem, we find that
$u_\infty$ belongs to $C^{2\alpha+\frac\epsilon2}(B_{\frac{|x_0|}4}(x_0))$. Then, $u_\infty$ is $C^{2\alpha+\frac\epsilon2}$ locally in $\R^N\setminus\{0\}$. Since $u_k$ is classical solution of (\ref{eq1.01}),  then it follows by Corollary 4.6 in \cite{CS1} that $u_\infty$ is a classical solution of (\ref{eq1.01}).
  \qquad$\Box$
\medskip

\noindent{\bf Proof of Theorem \ref{teo 2}.} The part $(i)$ follows by Proposition \ref{pr 3.3}. Now we prove
$(ii)$. By Proposition \ref{pr 3.4}, we see that $u_\infty$ is a classical solution of (\ref{eq1.01}).
It follows by uniqueness and rotation argument that
$u_k$ is radially symmetric.
The definition of $u_\infty$ and Proposition \ref{lm 13-09-1} yield
$u_\infty(x)=c_3|x|^{-\frac{2\alpha}{p-1}}$,  $x\in\R^N\setminus\{0\}$.\qquad$\Box$

\subsection{Properties of $u_{\infty,\Omega}$}

In this subsection, we make use of the properties of $\{u_k\}$ and $u_\infty$
 to estimate the weak solution of (\ref{eq1.001}) in general unbounded regular domain.

\smallskip
\noindent{\bf Proof of Theorem \ref{teo 3}.} First,
for any $k>0$, we use arguments in the proof of Lemma \ref{lm 3.1} to obtain that
 there exists a solution $u_{k,\Omega}$ of
(\ref{3.11}) replaced $B_R(0)$ by $\Omega$ 
such that
\begin{equation}\label{4.3.1}
u_k-m_{k,\Omega} \le u_{k,\Omega}\le u_k\quad {\rm in}\quad \Omega,
\end{equation}
where $m_{k,\Omega}=\sup_{x\in\Omega^c}u_k(x)$ and $u_k$ is given by Theorem \ref{cor 1}.
 Similarly to the proof of Proposition \ref{pr 3.1}, we obtain that $u_{k,\Omega}$ is  unique and then
\begin{equation}\label{5.2}
  u_{k,\Omega}=\lim_{R\to\infty} u_{k,\Omega\cap B_R(0)}.
\end{equation}
Since
$u_{k,\Omega\cap B_R(0)}$ is increasing respected to $k$ and $R$, the mapping $k\mapsto u_{k,\Omega}$ is increasing.

Next, by (\ref{5.2}) we have that $u_{\infty,\Omega}\ge u_{\infty,\Omega\cap B_R(0)}$ in $\R^N$ for any $R>0$.
For $p\in(0,\min\{1+\frac{2\alpha}N, \frac{N}{2\alpha}\})$,
using \cite[Theorem 1.1,Theorem 1.2]{CV3}, we have that
$\lim_{k\to\infty}u_{k,\Omega\cap B_R(0)}=\infty$ in $\Omega\cap B_R(0)$,
then $u_{\infty,\Omega}=\infty$ in $\Omega\cap B_R(0)$ for any $R>0$, which implies that
$$u_{\infty,\Omega}=\infty\quad{\rm in}\quad \Omega.$$

Finally, for $p\in(1+\frac{2\alpha}{N},\frac{N}{N-2\alpha})$, by (\ref{4.3.1}) we have that
$$u_\infty-m_{\infty,\Omega}\le u_{\infty,\Omega}\le u_\infty\quad {\rm in}\ \ \Omega,$$
where $m_{\infty,\Omega}=\sup_{x\in\Omega^c}u_\infty(x)\ge m_{k,\Omega}$, since $\{u_k\}$ are increasing.
Using arguments in the proof of Proposition \ref{pr 3.4}, we obtain that $u_{\infty,\Omega}$ is a classical solution of (\ref{eq1.02}). \qquad$\Box$


\begin{thebibliography}{99}
\bibitem {AH} D. R. Adams, L. I. Hedberg, Function spaces and potential theory, {\it Springer} (1996).

\bibitem {BLOP} D. Bartolucci, F. Leoni, L. Orsina and  A. C. Ponce, Semilinear equations with exponential nonlinearity and measure data,
{\it Ann.  I. H. Poincar\'{e}, 22(6)}, 799-815 (2005).

\bibitem {BB}
Ph. B\'{e}nilan and H. Brezis, Nonlinear problems related to the
Thomas-Fermi equation, {\it J. Evolution Eq., 3}, 673-770 (2003).

\bibitem {BBC} Ph. B\'{e}nilan, H. Brezis and M. Crandall, A semilinear elliptic equation in $L^1(\R^N )$, {\it Ann. Sc. Norm. Sup. Pisa Cl. Sci., 2}, 523-555 (1975).

\bibitem {B12} H. Brezis, Some variational problems of the Thomas-Fermi type.
Variational inequalities and complementarity problems, {\it Proc.
Internat. School, Erice, Wiley, Chichester}, 53-73 (1980).

\bibitem {BMP1}
H. Brezis, M. Marcus and  A. C. Ponce, A new concept of reduced measure for nonlinear elliptic equations.
{\it Comptes Rendus Mathematique, 339(3)}, 169-174 (2004).

\bibitem {BMP}
H. Brezis, M. Marcus and  A. C. Ponce, Nonlinear elliptic equations with measures revisited. Mathematical Aspects of Nonlinear Dispersive Equations, {\it Ann.  Math. Stud. 163}, 55-110 (2007).



\bibitem {BV}
H. Brezis and L. V\'{e}ron, Removable singularities of some nonlinear elliptic equations,
 {\it Arch. Rational Mech. Anal., 75}, 1-6 (1980).

\bibitem {CS2}
 L. Caffarelli and  L. Silvestre, An extension problem related to the fractional Laplacian.
 {\it Comm. Partial Differential Equations, 32(8)}, 1245-1260 (2007).

\bibitem {CS1}
 L. Caffarelli and  L. Silvestre, Regularity theory for fully non-linear integrodifferential equations,
 {\it Comm.   Pure   Appl. Math., 62}, 597-638 (2009).


 \bibitem {CS3} L. Caffarelli and L. Silvestre, Regularity results for nonlocal
equations by approximation, {\it Arch. Ration. Mech. Anal., 200(1),}
59-88, (2011).



\bibitem {CFQ} H. Chen, P. Felmer and A. Quaas, Large solution to elliptic  equations involving fractional Laplacian,
Accepted by {\it  Ann.  I. H. Poincar\'{e}.} (arXiv:1311.6044).


\bibitem {CKS} Z. Chen, P. Kim and R. Song,   Heat kernel estimates for the Dirichlet
fractional Laplacian. {\it J. Eur. Math. Soc., 12,}  1307-1329
(2010).

\bibitem {CS} Z. Chen, and R. Song, Estimates on Green functions and poisson kernels for symmetric stable process, {\it Math.
Ann., 312}, 465-501 (1998).

\bibitem {CT}Z. Chen and J. Tokle, Global heat kernel estimates  for fractional laplacians in unbounded open sets,
{\it Probab. Theory Related Field, 149,} 373-395 (2011).

\bibitem {CV1}  H. Chen and L. V\'{e}ron, Semilinear fractional elliptic equations with
gradient nonlinearity involving measures, {\it  J.  Funct. Anal., 266(8)}, 5467-5492 (2014).

\bibitem {CV2} H. Chen and L. V\'{e}ron, Semilinear fractional elliptic equations  involving measures,
{\it J. Diff. Eq.}, (arXiv:1305.0945).

\bibitem {CV3}  H. Chen and L. V\'{e}ron, Weakly and strongly singular solutions of semilinear fractional elliptic equations,
 {\it DOI 10.3233/ASY-141216,  (arXiv:1307.7023)}.



\bibitem {FQ2} P. Felmer and A. Quaas, Fundamental solutions and Liouville type theorems for nonlinear integral operators,
{\it Adv. Math., 226}, 2712-2738 (2011).


\bibitem {KPU} K. Karisen, F. Petitta  and S. Ulusoy, A duality approach to the fractional laplacian with measure data,
{\it Publ. Mat., 55}, 151-161 (2011).



\bibitem {MV3} M. Marcus  and L. V\'{e}ron, Removable singularities and boundary
traces, {\it J. Math. Pures Appl. 80}, 879-900 (2001).

\bibitem {MP} M. Marcus  and A. C. Ponce, Reduced limits for nonlinear equations with measures, {\it J. Funct. Anal., 258}, 2316-2372  (2010).



\bibitem {P} A. C. Ponce, Selected problems on elliptic equations involving measures.  ArXiv:1204.0668.


\bibitem {RS} X. Ros-Oton and J. Serra, The Dirichlet problem for the fractional
laplacian: regularity up to the boundary, {\it J. Math. Pures Appl., 101(3)},  275-302 (2014).

\bibitem {S} L. Silvestre, Regularity of the obstacle problem for a fractional power of the laplace operator,
{\it Comm. Pure Appl. Math., 60,} 67-112 (2007).

\bibitem {S1}
Y. Sire and E. Valdinoci, Fractional Laplacian phase transitions and
boundary reactions: a geometric inequality and a symmetry result,
{\it J. Funct. Anal., 256}, 1842-1864 (2009).


\bibitem {St} E. Stein, Singular integrals and differentiability properties of functions, {\it  Princeton University Press} (1970).



\bibitem {V1} L. V\'{e}ron, Weak and strong singularities of nonlinear elliptic equations, {\it Proc. Symp. Pure Math.,
45,} 477-495 (1986).

\bibitem {V0}  L. V\'{e}ron, Singular solutions of some nonlinear elliptic equations,
{\it Nonlinear Anal. T. M. $\&$ A., 5,} 225-242 (1981).

\bibitem {V}  L. V\'{e}ron, Elliptic equations involving Measures,
 Stationary Partial Differential equations,
{\it Vol. I, 593-712, Handb. Differ. Equ., North-Holland, Amsterdam}
(2004).




\end{thebibliography}
\end{document}